\documentclass[12pt, twoside]{article}
\usepackage{amsmath,amsthm,amssymb}
\usepackage{times}
\usepackage{enumerate}
\usepackage{extarrows}
\usepackage{color}
\usepackage{tabularx}
\usepackage{multirow}
\usepackage{graphicx}

\def\titlerunning#1{\gdef\titrun{#1}}
\makeatletter
\def\author#1{\gdef\autrun{\def\and{\unskip, }#1}\gdef\@author{#1}}
\def\address#1{{\def\and{\\\hspace*{18pt}}\renewcommand{\thefootnote}{}%
\footnote {#1}}%
\markboth{\autrun}{\titrun}}
\makeatother
\def\email#1{e-mail: #1}

\def\keywords#1{\par\medskip
\noindent\textbf{Keywords.} #1}

\newtheorem{theorem}{Theorem}[section]
\newtheorem{lemma}{Lemma}[section]

\newtheorem{remark}{Remark}[section]
\newtheorem{corollary}{Corollary}[section]

\newcommand{\Proof}{\begin{proof}}
\newcommand{\End}{\end{proof}}

\numberwithin{equation}{section}

\newcommand{\PreserveBackslash}[1]{\let\temp=\\#1\let\\=\temp}
\newcolumntype{C}[1]{>{\PreserveBackslash\centering}p{#1}}
\newcolumntype{R}[1]{>{\PreserveBackslash\raggedleft}p{#1}}
\newcolumntype{L}[1]{>{\PreserveBackslash\raggedright}p{#1}}

\newcolumntype{I}{!{\vrule width 1pt}}
\newlength\savedwidth

\begin{document}


\baselineskip=15pt


\titlerunning{ }

\title{A KAM theorem of symplectic algorithms for nearly integrabel Hamiltonian systems}

\author{Zaijiu Shang ,  Yang Xu}

\date{\today}

\maketitle

\address{Zaijiu Shang:  Institute of Academy of Mathematics and Systems Science, Chinese Academy of Sciences, Beijing 100080, China; \email{zaijiu@amss.ac.cn}
		\and Yang Xu: School of Mathematical Sciences, Fudan University, Shanghai 200433, China; \email{xuyang$\_$@fudan.edu.cn}}

\begin{abstract}
In this paper we prove a KAM-like theorem of symplectic algorithms for nearly integrable Hamiltonian systems which generalises the result of \cite{r1} and \cite{r6} for the case of integrable systems.
 
%
\end{abstract}

\tableofcontents

\section{Introduction}
\label{s1}
In \cite{r1} the author established a KAM theorem of symplectic algorithms for integrable Hamiltonian systems and left the question of generalization to nearly integrable systems open. In this paper we will solve it. The notations of \cite{r1} will be still used.

Consider the nearly integrable Hamiltonian system :
\begin{equation} \label{1}
    \dot{p} = - \frac{\partial H^\epsilon} {\partial q} (p,q) , \   \dot{q} =  \frac{\partial H^\epsilon}{\partial p} (p,q),  \   (p,q) \in D,
\end{equation}
where $H^\epsilon(p,q) = H^0 (p) + \epsilon H(p,q,\epsilon)$ is the Hamiltonian function,
$D \subseteq {\mathbb{R}}^{n} \times \mathbb{T}^n$ is a bounded  connected open domain, and the dot is the derivative with respect to time $t$.
For such nearly integrable Hamiltonian system, if the perturbation is small enough and the Hamiltonian function is sufficiently smooth, then there exist invariant tori, which is what we know. 
In this paper, we will use the method of discretization of symplectic algorithm to show that the system after the discretization of symplectic algorithm also has invariant tori.
Compared with the invariant tori of the nearly integrable Hamiltonian system, there is only slight deformation between them.
Without loss of generality, we can assume $ \lvert H  \rvert _ D \leq M_1$, where $\lvert \cdot \lvert_ D$ denotes the supremum norm on $D$, $M_1$ is a constant, and we can ignore the parameter expression of $\epsilon$. 

Denote $\omega$ as the derivative of  $H^0$, $\omega(\xi) =\frac{\partial H^0}{\partial \xi}(\xi)$, which is defined on $V_{\kappa} $, where $ V \subset \mathbb{R}^n $ is a bounded  connected domain,  $V_{\kappa} = V + \kappa := \bigcup \limits_{b \in V} \{ \xi \in \mathbb{C}^n : \lvert \xi - b \rvert_2 < \kappa \} \subseteq \mathbb{C}^n$, and $\kappa$ is a constant.

Use symplectic algorithm to discretize the Hamiltonian system \eqref{1}. By Lemma 3.1 and Lemma 3.3 in \cite{r1}, there is a function $P$, which depends on time step $t$, and it is well defined and real analytic on $ D_{\frac{r}{4},\frac{s}{4}}$, such that if time step $t$ is small enough, $0 <t < \delta$, where $\delta$ is a small enough constant, then the symplectic differential form of the system \eqref{1} after discretization by symplectic algorithm can be expressed as  $G_{H^ \epsilon}^t : (p,q) \to (\hat{p}, \hat{q}) $ 
\begin{equation} \label{2}
    G_{H^ \epsilon}^t :
    \begin{cases}
    \hat{p} = p 
    - t \epsilon \frac{\partial   S} {\partial q} (\hat{p}, q) 
    - t^{\alpha +1} \frac{ \partial{P}}{\partial q}(\hat{p}, q), \\
    \hat{q} = q + t \omega  
    + t \epsilon \frac{\partial   S} {\partial \hat{p}} (\hat{p}, q) 
    + t^{\alpha +1} \frac{ \partial{P}}{\partial \hat{p}}(\hat{p}, q),
    \end{cases}
\end{equation}
where $(\hat{p}, q) \in D_{\frac{r}{4},\frac{s}{4}}$, and $\alpha$ is a positive constant.
Here $S$ is a generating function, $S(p,q,t) = \sum \limits_{k=0}^\infty  (S)_k(p,q) t^k = H^1(p,q)  +  (S)_1(p,q) t +  (S)_2(p,q) t^2 + \dots$, and $D_{r,s}$ refers to
\begin{equation*}
    D_{r,s} = \{p : \inf \limits_{p^* \in V} \lvert p-p^* \rvert < r \} \times \{q : \Re q \in \mathbb{T}^n, \lvert \Im q \rvert < s \}
    \subset \mathbb{C}^n \times \mathbb{C}^n.
\end{equation*}
By Lemma 3.1 and Lemma 3.3 in \cite{r1},
$ \big \lvert \frac{\partial P}{\partial {\hat{p}}}  \big \rvert _ {\frac{r}{4},\frac{s}{4}} \leq M_2 $, $\big \lvert \frac{\partial P}{\partial q} \big \rvert _ {\frac{r}{4},\frac{s}{4}} \leq M_2 $, where $\lvert \cdot \lvert_ {\frac{r}{4},\frac{s}{4}}$ denotes the supremum norm on $D_{\frac{r}{4},\frac{s}{4}}$, and $M_2$ is a constant.
Similar to  \cite{r1},  fixed $(\hat{p_0}, q_0) \in D$, let $P(\hat{p_0}, q_0) =0$, for
$(\hat{p}, q) \in D_{\frac{r}{4},\frac{s}{4}}$, then 
we have $\lvert  P \rvert _ {\frac{r}{4},\frac{s}{4}} \leq 2 n M_2 l_* $, where $l_*$ is an upper bound of the length of the shortest curves from $(\hat{p_0}, q_0)$ to $(\hat{p}, q)$ in $D_{\frac{r}{4},\frac{s}{4}}$.
And, $\lvert  G_{H^ \epsilon}^t(p,q) -  g_{H^ \epsilon}^t(p,q) \rvert \leq M_3 t^{\alpha+1}$, where $g_{H^ \epsilon}^t$ is the phase flow of the nearly integrable Hamiltonian system \eqref{1}, and $M_3$ is a constant that does not depend on $t$.

For the nearly integrable Hamiltonian system \eqref{1}, there exists a generating function $\widetilde{S}$ such that its invariant tori can be found by the following system:
\begin{equation}  \label{xy3}
    \hat{p} = p 
    - t \epsilon \frac{\partial  \widetilde{S}} {\partial q} (\hat{p}, q) 
    , \ 
    \hat{q} = q + t \omega  
    + t \epsilon \frac{\partial  \widetilde{S}} {\partial \hat{p}} (\hat{p}, q) ,
\end{equation}
and for the time step $t$ in the common area, we have $S(p,q,t) - \widetilde{S}(p,q,t) \sim o(t^\alpha)$, let us denote it as $\lvert S - \widetilde{S} \rvert \leq M_4 t^\alpha$, where $M_4$ is a constant that does not depend on $t$.

In this paper, we assume that $\omega$ satisfies the R\"{u}ssmann's non-degeneracy condition.

\begin{remark} \label{re1}
$\omega$ satisfies the R\"{u}ssmann's non-degeneracy condition means that $\omega \in \{ \omega(\xi): \xi \in V \}$ does not lie in any hyperplane that passes through the origin, i.e. $\langle \omega , y \rangle \not\equiv 0 $ for all $y \in \mathbb{R}^n \setminus \{0\}$.
Xu, You and Qiu gave an expression of the R\"{u}ssmann's non-degeneracy condition for analytic cases(see Remark 3.1 in \cite{r2}),
that is, if $\omega$ satisfies the R\"{u}ssmann's non-degeneracy condition, then there exists an integer $\bar{n} >0$, such that
 \begin{equation} \label{4}
     rank \big \{\partial^i_{\xi}\omega(\xi) :  |i| \leq \bar{n} \big \} =n ,  \   \forall \xi \in V.
 \end{equation}
And, by R\"{u}ssmann (Lemma 18.2 in \cite{r5}), there are numbers $ \bar{n}(V) \in \mathbb{N}$  and $ \beta_0(V) > 0$ such that
\begin{equation}
    \min_{\xi \in V} \max_{0\leq v\leq \bar{n}} |D^v \langle k,  \omega(\xi) \rangle |\geq \beta_0,
\end{equation}
where $\langle k,  \omega \rangle = \sum  \limits_{j=1}^n k_j \omega_j$, $k \in \{c= (c_1, \dots, c_n) \in \mathbb{R}^n : \lvert c_j \rvert =1, j= 1, \dots, n \}$. We can take the smallest of such integer  $\bar{n} $.
Moreover, by R\"{u}ssmann (Theorem 18.4 in \cite{r5}), we have
\begin{equation}
    \min_{\xi \in V} \max_{0\leq v\leq \bar{n}} \big \lvert D^v {\lvert k \rvert_2}^{-2}  \lvert \langle k,  \omega(\xi) \rangle \rvert^2 \big \rvert \geq \beta_0,
\end{equation}
for $\forall k \in \mathbb{Z}^n \setminus \{0 \}$.
As $\lvert \cdot \rvert$ and $\lvert \cdot \rvert_2$ are equivalent, where $\lvert k \rvert = \sum  \limits_{j=1}^n \lvert k_j \rvert$, $\lvert k \rvert_2 = (\sum  \limits_{j=1}^n \lvert k_j \rvert^2)^{\frac{1}{2}}$, then there exist $\bar{n}$ and $\beta$ such that
\begin{equation}  \label{30}
    \min_{\xi \in V} \max_{0\leq v\leq \bar{n}} |D^v  \langle k,  \omega(\xi) \rangle |\geq \beta \lvert k \rvert,
\end{equation}
for $\forall k \in \mathbb{Z}^n \setminus \{0 \}$, where $\bar{n} = \bar{n}(\omega, V) \in \mathbb{N} $ is called the index of non-degeneracy of $\omega$ with respect to $V$, and $\beta = \beta(\omega, V) > 0$ is called the amount of non-degeneracy of $\omega$ with respect to $V$.
\end{remark}

Fixed $\tau \geq (n+2)(\bar{n} +1)$, $\gamma > 0$ and $0 < t < 1$, introduce the concept of Diophantine condition.
The Diophantine condition refers to that $\omega(\xi)$ satisfies 
\begin{equation} \label{5}
    \lvert e^{i \langle k, t \omega \rangle} - 1 \rvert \geq \frac{t \gamma}{{\lvert k \rvert}^\tau} , \ \forall \  k = (k_1, k_2, \dots , k_n ) \in \mathbb{Z}^n \setminus \{0 \},
\end{equation}
where $\xi \in V_\kappa$, $\omega = (\omega_1, \omega_2, \dots , \omega_n) $ , $\langle k,  \omega \rangle = \sum  \limits_{j=1}^n k_j \omega_j$ , $\lvert k \rvert = \sum  \limits_{j=1}^n \lvert k_j \rvert$.  $\lvert \cdot \lvert_ {r,s}$ denotes the supremum norm on $D_{r,s} $.

\section{Main theorem}
\label{s2}

\begin{theorem} \label{t1}
For the nearly integrable Hamiltonian system \eqref{1} , $H $ is real analytic on  $D_{r,s} \times V_\kappa $, and $H^0 $is analytic on $V_\kappa $.
If time step  $t$ and disturbance parameter  $\epsilon$ are small enough, $\omega$ satisfies the R\"{u}ssmann's non-degeneracy condition, 
then the generating function representation of symplectic algorithm \eqref{2} has invariant tori on a set of large  measures. 
The generating function representation of phase flow in 
nearly integrable Hamiltonian system \eqref{xy3} also has invariant tori on a set of large  measures, and there is only a slight deformation in the common area. That is to say, for system \eqref{2}, there is a non-empty Cantor set $V_{\epsilon,t} \subseteq V$ and a Whitney smooth symplectic mapping  
$\Phi _{\epsilon,t} : V_{\epsilon,t} \times \mathbb{T}^n \to \mathbb{R}^n \times \mathbb{T}^n $ that makes the following true.

\begin{enumerate}
    \item[(\romannumeral1)] 
    $\Phi _{\epsilon,t}$ is a symplectic conjugation between $G_{H^ \epsilon}^t$ and $R_{\epsilon,t}$, i.e.,
    $\Phi _{\epsilon,t}^{-1} \circ G_{H^ \epsilon}^t \circ \Phi _{\epsilon,t} = R_{\epsilon,t}$, where $R_{\epsilon,t}$ is a rotation on $V_{\epsilon,t} \times \mathbb{T}^n$ with frequency $t \omega_{\epsilon,t} $, i.e., $R_{\epsilon,t}(p,q) = (p, q + t \omega_{\epsilon,t} )$.
    \item[(\romannumeral2)] 
    $V_{\epsilon,t}$ is a set of positive measure if $\gamma$ is small enough. And for $\gamma \to 0$, we have $\lvert V \setminus V_{\epsilon,t} \rvert \to 0 $.
    To be specific, 
    \begin{equation} \label{6}
        \lvert V \setminus V_{\epsilon,t} \rvert \leq c_7   \gamma  d^{n}((n+1)^{-\frac{1}{2}}+2d+ \kappa^{-1}d),
    \end{equation}
     where $c_7 = 12 \bar{c} (2\pi e)^\frac{n+1}{2} (\bar n +2)^{\bar n +3} [(\bar n  +2)!]^{-1} \hat{K} \beta^{-\frac{\bar{n}+2}{\bar{n } +1}} \lvert  \omega \rvert^{\bar{ n}  +2}_{\mathcal{B}} $, $\hat{K} = \lvert \omega \rvert_{V} + c_4^{'}( M_1  +  2 n M_2 l_* )   +1$, $\mathcal{B} = (V \times (0,1) + \kappa) \cap \mathbb{R}^{n+1}$, $d$ is the diameter of $V \times (0,1)$, $\bar{c}$ is a constant depending only on $n$, $\bar{n}$ and $\tau$, and $c_4^{'}$ is  a constant depending only on $\tau, n, \bar{n}, \gamma,  \rho_0, \sigma_0$ and $K_0$.
     \item[(\romannumeral3)] 
         Let the time step be $t_1$ and $t_2$ respectively, and comppare the two systems, then for $\xi \in V_{\epsilon,t_1} \cap V_{\epsilon,t_2}$ we have
    \begin{equation*}
    \lvert \Phi_{\epsilon,t_1} - \Phi_{\epsilon,t_2} \rvert \leq  2 c_2^{'}  n M_2 l_* (t_1^{\alpha}-t_2^{\alpha})  , \ \lvert \omega_{\epsilon,t_1} - \omega_{\epsilon,t_2} \rvert \leq  2 c_4^{'} n M_2 l_* (t_1^{\alpha}-t_2^{\alpha}),
\end{equation*}
   where $M_2$ and $l_*$ are constants, and $c_2^{'}$ and $c_4^{'}$ are constants depending only on $\tau, n, \bar{n},  \rho_0, \sigma_0, \gamma_0 $ and $K_0$. 
$ V_{\epsilon,t_1} \cap V_{\epsilon,t_2} $ is also a set of positive measure if $\gamma$ is small enough. And for $\gamma \to 0$, we also have $\lvert V \setminus (V_{\epsilon,t_1} \cap V_{\epsilon,t_2}) \rvert \to 0 $.
\end{enumerate}
\end{theorem}

\section{Proof of Theorem \ref{t1}}
\label{s3}

\subsection{One-step Iterative Analysis}
Firstly, we focus on one-step iteration.
Set
\begin{equation*}
    V_{\gamma,t,v} = \big \{\xi \in V : \lvert e^{i \langle k, t \omega_v(\xi) \rangle} - 1 \rvert \geq \frac{t \gamma_v}{{\lvert k \rvert}^\tau} , \ \forall \  k \in \mathbb{Z}^n \setminus \{0 \} \big \},
\end{equation*}
that is, for $\xi \in V_{\gamma,t,v}$ , we have $\omega_v(\xi)$ satisfies the Diophantine condition \eqref{5}.
It is a Cantor set, about which measure estimates are given later. Whitney analytical extension was carried out for $\omega_v(\xi)$ on $ V_{\gamma,t,v}$ and denoted as $V_{\kappa^{'}_v} = V_{\gamma,t,v} + \kappa^{'}_v$.

Denote the one-step transforming of the symplectic difference scheme as
$G_v : (p,q) \to (\hat{p}, \hat{q})$ with
\begin{equation} \label{7}
    G_v : 
    \begin{cases}
        \hat{p} = p 
    - \partial_2 \big( t \epsilon S_v(\hat{p}, q) 
    + t^{\alpha +1} P_v(\hat{p}, q)\big) ,\\
    \hat{q} = q + t \omega_v  
    + \partial_1 \big(t \epsilon S_v(\hat{p}, q) 
    + t^{\alpha +1} P_v(\hat{p}, q)\big),
    \end{cases}
\end{equation}
here the subscript $v$ represents the $v$-th step.
Assume the one-step symplectic mapping  can be expressed through the generating function $\psi_v$ as $\Psi_v : (I,\theta) \to (p,q)$ with
\begin{equation} \label{8}
    \Psi_v : 
    \begin{cases}
     I = p  - \partial_2 \psi_v (I,q) ,\\
     \theta = q  + \partial_1 \psi_v, (I,q),
    \end{cases}
\end{equation}
so we have $(I,\theta) \stackrel{\Psi_v}{\longrightarrow} (p,q) \stackrel{G_v}{\longrightarrow} (\hat{p},\hat{q}) \stackrel{\Psi_v^{-1}}{\longrightarrow} (\hat{I},\hat{\theta})$, and $\Psi_v^{-1} \circ G_v \circ \Psi_v $ can be expressed as follows:
\begin{align*}
    \hat{I} &= \hat{p} - \partial_2 \psi_v (\hat{I},\hat{q}) \\
    &= p - \partial_2 \big(t \epsilon S_v(\hat{p}, q) + t^{\alpha +1} P_v(\hat{p}, q)\big) - \partial_2 \psi_v (\hat{I},\hat{q})\\
    &= I - \partial_2 \big(t \epsilon S_v(\hat{p}, q) 
    + t^{\alpha +1} P_v(\hat{p}, q)\big) + \partial_2 \psi_v (I,q) - \partial_2 \psi_v (\hat{I},\hat{q}),\\
    \hat{\theta} =& \ \hat{q} + \partial_1 \psi_v (\hat{I},\hat{q}) \\
    =& \ q + t \omega_v + \partial_1 \big(t \epsilon S_v(\hat{p}, q) + t^{\alpha +1} P_v(\hat{p}, q)\big) + \partial_1 \psi_v (\hat{I},\hat{q})\\
    =& \ \theta + t \omega_v  + \partial_1 \big(t \epsilon S_v(\hat{p}, q) 
    + t^{\alpha +1} P_v(\hat{p}, q)\big) - \partial_1 \psi_v (I,q) + \partial_1 \psi_v (\hat{I},\hat{q}).
\end{align*}
Note that for $0< t < \delta$, we have
\begin{equation*}
    \big \lvert \partial_1 \big(t \epsilon S_0(\hat{p}, q) + t^{\alpha +1} P_0(\hat{p}, q)\big) \big \rvert_{\frac{r}{4},\frac{s}{4}} \leq \frac{ M_1 t\epsilon +  2 n M_2 l_* t^{\alpha +1}}{r / 4},
\end{equation*}
\begin{equation*}
    \big \lvert \partial_2 \big(t \epsilon S_0(\hat{p}, q) + t^{\alpha +1} P_0(\hat{p}, q)]\big) \big \rvert_{\frac{r}{4},\frac{s}{4}} \leq \frac{ M_1 t\epsilon +  2 n M_2 l_* t^{\alpha +1}}{s / 4}.
\end{equation*}
Thus, for $\epsilon$ and $t$ small enough such that $ \frac{ M_1 \epsilon +  2 n M_2 l_* t^{\alpha}}{r / 4} , \frac{ M_1 \epsilon +  2 n M_2 l_* t^{\alpha}}{s / 4} \leq E_0$, then, if $r_0 \leq \frac{r}{4}$, $s_0 \leq \frac{s}{4}$,  we can denote $\big \lvert \partial_1 (t \epsilon S_0 + t^{\alpha +1} P_0) \big \rvert _{r_0,s_0} \leq t E_0$,
$\big \lvert \partial_2 (t \epsilon S_0 + t^{\alpha +1} P_0)  \big \lvert _{r_0,s_0} \leq t E_0$ , here $E_0$ is also small enough as $\epsilon$ and $t$ are small enough.

Similarly, we can assume that
$\big \lvert \partial_1 (t \epsilon S_v + t^{\alpha +1} P_v) \big \rvert _{r_v,s_v} \leq t E_v$  as well as
$\big \lvert \partial_2 (t \epsilon S_v + t^{\alpha +1} P_v)  \big \lvert _{r_v,s_v}  \leq t E_v$, where $E_v$ is also a small quantity.
Then, as $\psi_v$ and $(t \epsilon S_v + t^{\alpha +1} P_v)$ are small quantities, we can get the following analysis:
\begin{align*}
    \partial_1 \psi_v (I,q) 
    = \ & \partial_1 \psi_v \big(I,\theta - \partial_1 \psi_v (I,q)\big)
    =  \partial_1 \psi_v (I,\theta) + o(E_v) , \\
    \partial_1 \psi_v (\hat{I},\hat{q})
    = \ &
    \partial_1 \psi_v \big(I - \partial_2 (t \epsilon S_v 
    + t^{\alpha +1} P_v)(\hat{p}, q) + \partial_2 \psi_v (I,q) - \partial_2 \psi_v (\hat{I},\hat{q}), \\
    & \theta + t \omega_v  + \partial_1 (t \epsilon S_v 
    + t^{\alpha +1} P_v)(\hat{p}, q) - \partial_1 \psi_v (I,q) \big) \\
    = \ & \partial_1 \psi_v \big(I, \theta + t \omega_v (\xi) \big) + o(E_v) ,\\
    \partial_1 (t \epsilon S_v 
    + t^{\alpha +1}& P_v) (\hat{p}, q) =  \partial_1 (t \epsilon S_v 
    + t^{\alpha +1} P_v) \big( I - \partial_2 (t \epsilon S_v + t^{\alpha +1} P_v)(\hat{p}, q)  \\
    & \ + \partial_2 \psi_v (I,q), \theta - \partial_1 \psi_v (I,q)\big) \\
    = \ & \partial_1 (t \epsilon S_v 
    + t^{\alpha +1} P_v)( I, \theta) + o(E_v).
\end{align*}
Similarly, we have
\begin{align*}
    \partial_2 \psi_v (I,q) 
    = \ & \partial_2 \psi_v \big(I,\theta - \partial_1 \psi_v (I,q) \big)
    =  \partial_2 \psi_v (I,\theta) + o(E_v) ,\\
    \partial_2 \psi_v (\hat{I},\hat{q}) 
    = \ &
    \partial_2 \psi_v \big(I - \partial_2 (t \epsilon S_v 
    + t^{\alpha +1} P_v)(\hat{p}, q) + \partial_2 \psi_v (I,q) - \partial_2 \psi_v (\hat{I},\hat{q}), \\
    & \theta + t \omega_v  + \partial_1 (t \epsilon S_v 
    + t^{\alpha +1} P_v)(\hat{p}, q) - \partial_1 \psi_v (I,q) \big) \\
    = \ & \partial_2 \psi_v \big(I, \theta + t \omega_v (\xi) \big) + o(E_v), \\
    \partial_2 (t \epsilon S_v 
    + t^{\alpha +1}& P_v) (\hat{p}, q) =  \partial_2 (t \epsilon S_v 
    + t^{\alpha +1} P_v) \big( I - \partial_2 (t \epsilon S_v + t^{\alpha +1} P_v)(\hat{p}, q)  \\
    & \ + \partial_2 \psi_v (I,q), \theta - \partial_1 \psi_v (I,q) \big) \\
    = \ & \partial_2 (t \epsilon S_v 
    + t^{\alpha +1} P_v)( I, \theta) + o(E_v).
\end{align*}
So, we have $\Psi_v^{-1} \circ G_v \circ \Psi_v $ can be expressed as follows:
\begin{align*}
    \hat{I} = & \ I - \partial_2 \big(t \epsilon S_v(I,\theta)  + t^{\alpha +1} P_v(I,\theta) \big) + \partial_2 \psi_v (I,\theta) \\
     & - \partial_2 \psi_v \big(I,\theta + t \omega_v \big)+ o(E_v), \\
     \hat{\theta} = & \ \theta + t \omega_v(\xi) + \partial_1 \big(t \epsilon S_v(I,\theta)  + t^{\alpha +1} P_v(I,\theta) \big) - \partial_1 \psi_v (I,\theta) \\
     & + \partial_1 \psi_v \big(I,\theta + t \omega_v \big)+ o(E_v).
\end{align*}

Looking at this expression, we can truncate $t \epsilon S_v + t^{\alpha +1} P_v$ by zero degree term and one degree term with respect to the action variables to get $(t \epsilon S_v + t^{\alpha +1} P_v)^*$. 
For general function $f(I,\theta)$, entry notation $[f]$ and Fourier coefficients $f_k(I)$,
\begin{align*}
    & f(I,\theta)= \sum \limits _{k \in \mathbb{Z}^n} f_k(I) e^{i \langle k, \theta \rangle} = \widetilde{f(I,\theta)} + [f], \\
   & \widetilde{f(I,\theta)} = \sum \limits _{k \in \mathbb{Z}^n \setminus \{0 \}} f_k(I) e^{i \langle k, \theta \rangle}.
\end{align*}
To get $\psi_v$, we consider the following equation:
\begin{equation*}
     \psi_v \big(I,\theta + t \omega_v \big) -  \psi_v (I,\theta) +  \widetilde{\big(t \epsilon S_v(I,\theta)  + t^{\alpha +1} P_v(I,\theta) \big)^*} = 0,
\end{equation*}
where
\begin{align*}
     & \  \widetilde{ \big(t \epsilon S_v(I,\theta)  + t^{\alpha +1} P_v(I,\theta) \big)^*} \\
     =  & \ \big(t \epsilon S_v(I,\theta)  + t^{\alpha +1} P_v(I,\theta) \big)^* - \big [   \big(t \epsilon S_v(I,\theta)  + t^{\alpha +1} P_v(I,\theta) \big)^* \big ] \\
     =  & \sum \limits_{k \in \mathbb{Z}^n \setminus \{0 \}} \big(  (t \epsilon S_v  + t^{\alpha +1} P_v )^* \big)_k (I) e^{i \langle k, \theta \rangle}.
\end{align*}
Here $\big [  (t \epsilon S_v  + t^{\alpha +1} P_v )^* \big ]$  is the mean values of $  (t \epsilon S_v  + t^{\alpha +1} P_v )^* $   with respect to the angle variables on $\mathbb{T}^n$, i.e.
\begin{align*}
    \big [   (t \epsilon S_v  + t^{\alpha +1} P_v )^* \big ] & = \frac{1}{(2\pi)^n} \int_{\mathbb{T}^n}  \big(t \epsilon S_v(I,\theta)  + t^{\alpha +1} P_v(I,\theta) \big)^*  d \theta .
\end{align*}

Define $\omega_{v+1} = \omega_v +  \big [ \partial_1  (\epsilon S_v  + t^{\alpha} P_v )^* \big (\xi)]$. By the way, $[ \partial_2  (\epsilon S_v  + t^{\alpha} P_v )^* \big ] = 0$.
Moreover, we truncate the Fourier series expansion of $ (t \epsilon S_v  + t^{\alpha +1} P_v )^*$ with respect to angle variables by order $K_v$.
And then we get the final homological equation:
\begin{equation} \label{50}
     \psi_v (I,\theta + t \omega_v ) -  \psi_v (I,\theta) + T_{K_v} \widetilde{\big(t \epsilon S_v(I,\theta)  + t^{\alpha +1} P_v(I,\theta) \big)^*} = 0,
\end{equation}
where
\begin{align*}
    T_{K_v} \widetilde{\big(t \epsilon S_v(I,\theta)  + t^{\alpha +1} P_v(I,\theta) \big)^*} = \sum \limits_{0< \lvert k \rvert \leq K_v}  \big(  (t \epsilon S_v  + t^{\alpha +1} P_v )^* \big)_k e^{i \langle k, \theta \rangle}.
\end{align*}

To solve $ \psi_v$ in the homological equation \eqref{50}, we can expand $ \psi_v (I,\theta + t \omega_v )$, $ \psi_v (I,\theta)$ and $ \big(t \epsilon S_v(I,\theta)  + t^{\alpha +1} P_v(I,\theta) \big)^* $ by Fourier series, and then compare their coefficients. Their Fourier series expansion are shown below:
\begin{align*}
      \psi_v (I,\theta + t \omega_v )  = & [ \psi_v](I) + \sum \limits_{k \in \mathbb{Z}^n \setminus \{0 \}} \big( \psi_v \big)_k e^{i \langle k, \theta + t \omega_v \rangle},\\
      \psi_v (I,\theta )  =  [ \psi_v](I) & + \sum \limits_{k \in \mathbb{Z}^n \setminus \{0 \}} \big( \psi_v \big)_k e^{i \langle k, \theta \rangle},\\ 
      \big(t \epsilon S_v(I,\theta)  + t^{\alpha +1} P_v & (I,\theta) \big)^* =   \big [   (t \epsilon S_v  + t^{\alpha +1} P_v )^* \big ](I) \\
     & + \sum \limits_{k \in \mathbb{Z}^n \setminus \{0 \}} \big(  (t \epsilon S_v  + t^{\alpha +1} P_v )^* \big )_k e^{i \langle k, \theta \rangle}.
\end{align*}
Compare their coefficients of $e^{i \langle k, \theta \rangle}$, then we have $\big( \psi_v \big)_k = \frac{- \big(  (t \epsilon S_v  + t^{\alpha +1} P_v )^* \big )_k}{e^{i \langle k, t \omega_v \rangle} -1}$, so we get $ \psi_v$ :
\begin{equation} \label{9}
     \psi_v(I,\theta) = - \sum \limits_{0< \lvert k \rvert \leq K_v} \frac{ \big(  (t \epsilon S_v  + t^{\alpha +1} P_v )^* \big )_k e^{i \langle k, \theta \rangle}} {e^{i \langle k, t \omega_v \rangle} -1}.
\end{equation}

\begin{lemma} \label{l5}
    (see Lemma A.1 in \cite{r4}) If $f \in A^s$, then $f= \sum \limits _{k} f_k e^{i \langle k, \theta \rangle}$ with
   \begin{equation*}
       \lvert f_k \rvert \leq \lvert f \rvert_s e^{- \lvert k \rvert s} ,  \     k \in \mathbb{Z}^n
   \end{equation*}
  where
  $A^s$ is the set of all real analytic functions on $\{\theta : \lvert \Im (\theta) \rvert < s \}
    \subset \mathbb{C}^n$ with sup-norm $\lvert \cdot \rvert_s$, $\lvert k \rvert = \lvert k_1 \rvert + \dots + \lvert k_n \rvert$.
\end{lemma}

By Lemma \ref{l5}, $\lvert \big(\partial_1  (t \epsilon S_v  + t^{\alpha +1} P_v )^* \big )_k \rvert \leq \lvert \partial_1  (t \epsilon S_v  + t^{\alpha +1} P_v )^* \rvert_{s_v} e^{- \lvert k \rvert s_v}$.  $\omega_v$ satisfies the Diophantine condition \eqref{5} on $V_{\gamma,t,v}$, and we have
\begin{align*}
    \lvert  \partial_1 \psi_v \rvert _{s_v - \sigma_v} 
    & \leq \sum \limits_{0< \lvert k \rvert \leq K_v}
    \frac{ \lvert \partial_1 \big( (t \epsilon S_v  + t^{\alpha +1} P_v )^* \big )_k e^{i \langle k, \theta \rangle} \rvert_{{s_v - \sigma_v} }}{\lvert e^{i \langle k, t \omega_v \rangle} -1 \rvert} \\
    & \leq \sum \limits_{0< \lvert k \rvert \leq K_v}
    \frac{\lvert \big(\partial_1  (t \epsilon S_v  + t^{\alpha +1} P_v )^* \big )_k \rvert e^{ \lvert k \rvert (s_v - \sigma_v)}}{\lvert e^{i \langle k, t \omega_v \rangle} -1 \rvert} \\
    & \leq \frac{1}{t \gamma_v} \sum \limits_{0< \lvert k \rvert \leq K_v} \lvert \big(\partial_1  (t \epsilon S_v  + t^{\alpha +1} P_v )^* \big )_k \rvert e^{ \lvert k \rvert (s_v - \sigma_v)} \lvert k \rvert ^\tau \\
    & \leq \frac{ \lvert \partial_1  (t \epsilon S_v  + t^{\alpha +1} P_v )^* \rvert_{s_v}}{t \gamma_v} \sum \limits_{0< \lvert k \rvert \leq K_v} e^{- \lvert k \rvert   \sigma_v} \lvert k \rvert ^\tau \\
    & \leq \frac{c_1}{t \gamma_v \sigma_v ^{\tau +n}} \lvert \partial_1  (t \epsilon S_v  + t^{\alpha +1} P_v )^* \rvert_{s_v} \\
    & \leq \frac{2 c_1}{t \gamma_v \sigma_v ^{\tau +n}} \lvert \partial_1  (t \epsilon S_v  + t^{\alpha +1} P_v ) \rvert_{s_v} ,
\end{align*}
where $c_1$ is a constant only depends on $\tau$ and $n$. Similarly,  we have
\begin{align*}
    \lvert  \partial_2 \psi_v \rvert _{s_v - \sigma_v} 
    & \leq \sum \limits_{0< \lvert k \rvert \leq K_v}
    \frac{ \lvert \partial_2 \big( ( (t \epsilon S_v  + t^{\alpha +1} P_v )^*  )_k e^{i \langle k, \theta \rangle  }\big ) \rvert_{{s_v - \sigma_v} }}{\lvert e^{i \langle k, t \omega_v \rangle} -1 \rvert} \\
    & \leq \sum \limits_{0< \lvert k \rvert \leq K_v}
    \frac{ \lvert  k i  \big( (t \epsilon S_v  + t^{\alpha +1} P_v )^* \big )_k e^{i \langle k, \theta \rangle} \rvert_{{s_v - \sigma_v} } }{\lvert e^{i \langle k, t \omega_v \rangle} -1 \rvert  } \\
    & \leq \sum \limits_{0< \lvert k \rvert \leq K_v}
    \frac{ \lvert k i \big(  (t \epsilon S_v  + t^{\alpha +1} P_v )^* \big )_k \rvert e^{ \lvert k \rvert (s_v - \sigma_v)}}{\lvert e^{i \langle k, t \omega_v \rangle} -1 \rvert} \\
    & \leq \frac{1}{t \gamma_v} \sum \limits_{0< \lvert k \rvert \leq K_v} \lvert k i  \big(  (t \epsilon S_v  + t^{\alpha +1} P_v )^* \big )_k \rvert e^{ \lvert k \rvert (s_v - \sigma_v)} \lvert k \rvert ^{\tau} \\
     & \leq \frac{ \lvert \partial_2 (t \epsilon S_v  + t^{\alpha +1} P_v )^* \rvert_{s_v}}{t \gamma_v} \sum \limits_{0< \lvert k \rvert \leq K_v} e^{- \lvert k \rvert   \sigma_v} \lvert k \rvert ^{\tau}
\\
    & \leq \frac{c_1}{t \gamma_v \sigma_v ^{\tau +n}} \lvert \partial_2  (t \epsilon S_v  + t^{\alpha +1} P_v )^* \rvert_{s_v} \\
    & \leq \frac{2 c_1}{t \gamma_v \sigma_v ^{\tau +n}} \lvert \partial_2  (t \epsilon S_v  + t^{\alpha +1} P_v ) \rvert_{s_v}.
\end{align*}

Now we give the immature parameter settings, and the final complete parameter settings will come later.
\begin{equation*}
    \begin{array}{ccc}
    \tau_v = \frac{1^{-2}+\dots +v^{-2}}{2\sum \limits_{v=1}^\infty v^{-2}},
    & s_v = \frac{1}{4}(1-\tau_v) s_0 , \
    & \sigma_v = \frac{1}{4} (s_v - s_{v+1}), \\
    K_{v+1} = 4 K_v ,
    & r_{v+1}=\eta_v r_v ,  \
    & \rho_v = \frac{1}{4} (r_v - r_{v+1}), \\
     \gamma_v = \frac{\gamma^{\bar{n}+1}}{2^{(\bar{n}+1)v}},
    & F_v = \frac{E_v}{\gamma_v \sigma_v^{\tau+n+1}\rho_v} ,
    & F_{v+1} = F_v^{\frac{v+1}{v}}, \\
     \eta_v = \frac{1}{12 c_3} F_v^{\frac{1}{v}}.
\end{array}
\end{equation*}
Note: $c_3$  is a constant. In fact, $\eta_v$ is also a constant, as will be explained later.

Denote $\Phi_v = \Psi_0 \circ \Psi_1 \circ \dots \circ \Psi_{v-1}$, and let $C_v= R_v^{-1} \circ G_v$, where $R_v$ is a rotation with frequency $t \omega_v $, i.e., $R_v(I,\theta) = (I, \theta + t \omega_v )$.
Focus on one-step iteration, we have the following iteration lemma.

\begin{lemma}{(Iteration lemma)} \label{l1}

   Assume that $\big \lvert \partial_1 (t \epsilon S_v + t^{\alpha +1} P_v) \big \rvert _{r_v,s_v},  \big \lvert \partial_2 (t \epsilon S_v + t^{\alpha +1} P_v) \big \rvert _{r_v,s_v} \leq t E_v$, $t \omega_v$ satisfies the Diophantine condition \eqref{5} with $t \gamma_v$, and other parameters are the same as above. Take the appropriate initial values, then for $v \geq 3$, there exist a real analytic symplectic mapping $\Phi_v$ which is defined on $ D_{\eta_v r_v, s_v - 4 \sigma_v} \times V_{\kappa^{'}_v} $ such that $\Phi_v^{-1} \circ G_{H^\epsilon}^t \circ \Phi_v : (p,q) \to (\hat{p}, \hat{q})$
     \begin{equation*} 
    \begin{cases}
        \hat{p} = p 
    - \partial_2 \big( t \epsilon S_{v+1}(\hat{p}, q) 
    + t^{\alpha +1} P_{v+1}(\hat{p}, q)\big) ,\\
    \hat{q} = q + t \omega_{v+1}  
    + \partial_1 \big(t \epsilon S_{v+1}(\hat{p}, q) 
    + t^{\alpha +1} P_{v+1}(\hat{p}, q)\big),
    \end{cases}
  \end{equation*}
 where $\big \lvert \partial_1 (t \epsilon S_{v+1} + t^{\alpha +1} P_{v+1}) \big \rvert _{r_{v+1},s_{v+1}},  \big \lvert \partial_2 (t \epsilon S_{v+1} + t^{\alpha +1} P_{v+1}) \big \rvert _{r_{v+1},s_{v+1}} \leq t E_{v+1}$, and the following results are true.
 \begin{enumerate}
    \item[(\romannumeral1)]  \label{l11}
        $\Phi_v$ is well defined and real analytic on $D_{r_{v-1} - 3 \rho_{v-1}, s_{v-1} - 3 \sigma_{v-1}} \times V_{\kappa^{'}_v}$, and it has the following estimate:
        \begin{equation} \label{12}
           \lvert \Phi_v - \Phi_{v-1} \rvert \leq \frac{4 c_1  E_{v-1}}{\gamma_{v-1} \sigma_{v-1}^{\tau+n}},
        \end{equation}
        where $c_1$ is a constant only depends on $\tau$ and $n$.
    \item[(\romannumeral2)]  \label{l12}
        $C_v = R_v^{-1} \circ \Phi_v^{-1} \circ G_{H^\epsilon}^t \circ \Phi_v $ is well defined and real analytic on $D_{r_v, s_v} \times V_{\kappa^{'}_v} $, and it has the following estimate:
        \begin{equation} \label{13}
            \lvert C_v - i d \rvert _{r_v, s_v} \leq t E_v.
        \end{equation}
    \item[(\romannumeral3)]  \label{l13}
        For $\omega = \omega_0$, on $ V_{\kappa^{'}_v}$ we have
         \begin{equation} \label{14}
            \lvert \omega_v - \omega_{v-1} \rvert \leq 2 E_v.
        \end{equation}
 \end{enumerate}
\end{lemma}

\subsection{Proof of the Iteration Lemma}
Firstly, $C_v= R_v^{-1} \circ G_v$ :
\begin{equation*}
    \begin{cases}
     \hat{p} = p 
    - \partial_2 \big( t \epsilon S_v(\hat{p}, q) 
    + t^{\alpha +1} P_v(\hat{p}, q)\big), \\
    \hat{q} = q   
    + \partial_1 \big(t \epsilon S_v(\hat{p}, q) 
    + t^{\alpha +1} P_v(\hat{p}, q)\big).
    \end{cases}
\end{equation*}
In order to express $C_v$ explicitly, firstly, let us solve the first equation:
\begin{equation*}
    \hat{p} = p 
    - \partial_2 \big( t \epsilon S_v(\hat{p}, q) 
    + t^{\alpha +1} P_v(\hat{p}, q)\big).
\end{equation*}
Let $X_1$ be the set of all real analytic functions $p_*(p,q)$ on $D_{r_v, s_v} $, where $p_*(p,q)$ is periodic $2 \pi$ with respect to angle variables, and $ \lvert p_*(p,q) \rvert _{r_v, s_v} \leq \rho_v$.

Consider mapping:
\begin{equation*}
    f_1: p_* \to f_1(p_*) = - \partial_2 \big( t \epsilon S_v(p+ p_*, q) + t^{\alpha +1} P_v(p+p_*, q)\big).
\end{equation*}
We know that $f_1$ is well defined and real analytic on $X_1$, and $f_1 :X_1 \to X_1$. By Cauchy estimates, for $\forall \ p_* \in X_1$, $(p,q) \in D_{r_v, s_v}$we have
\begin{equation*}
    \lvert f_1(p_*) \rvert \leq \lvert \partial_2 \big( t \epsilon S_v(p+p_*, q) + t^{\alpha +1} P_v(p+p_*, q)\big) \rvert \leq \rho_v.
\end{equation*}
For $\forall \ p_*^1, p_*^2 \in X_1$, $(p,q) \in D_{r_v, s_v}$,  we have
\begin{align*}
    \lvert f_1(p_*^1) - f_1(p_*^2) \rvert  = &   \ \lvert \partial_2 \big( t \epsilon S_v(p+p_*^1, q) + t^{\alpha +1} P_v(p+p_*^1, q)\big) \\
    & -   \partial_2 \big( t \epsilon S_v(p+p_*^2, q) + t^{\alpha +1} P_v(p+p_*^2, q)\big) \rvert \\
     \leq &  \ n \lvert \partial_1 \partial_2 \big( t \epsilon S_v(p+\bar{p}, q) + t^{\alpha +1} P_v(p+\bar{p}, q)\big) \rvert \cdot \lvert p_*^1 - p_*^2 \rvert \\
     \leq & \ \frac{n t E_v}{\rho_v} \cdot \lvert p_*^1 - p_*^2 \rvert \\
     \leq & \ L_1 \lvert p_*^1 - p_*^2 \rvert,
\end{align*}
if $\frac{n t E_0}{\rho_0} <1 $, where $0<L_1<1, \bar{p} \in X_1$.

So, $f_1: X_1 \to X_1$ is a compressed mapping, by the fixed point theorem, there is a unique fixed point $p_{**} \in X_1$, such that  $C_v$ can be expressed explicitly as:
\begin{equation*}
    \begin{cases}
     \hat{p} = p + p_{**}, \\
    \hat{q} = q   
    + \partial_1 \big(t \epsilon S_v(p + p_{**}, q)
    + t^{\alpha +1} P_v(p + p_{**}, q)\big),
    \end{cases}
\end{equation*}
so we get the following corollary.
\begin{corollary}
   $\partial_1 (t \epsilon S_v 
    + t^{\alpha +1} P_v)(p+ p_{**},q,\xi)$ , $\partial_2 (t \epsilon S_v + t^{\alpha +1} P_v)(p+ p_{**},q,\xi)$ are well defined and real analytic on $D_1= D_{r_v -  \rho_v, s_v -  \sigma_v} \times V_{\kappa^{'}_v} $, and they are periodic $2 \pi$ with respect to angle variables. Moreover, there are the following estimates on $D_1$:
    \begin{equation*}
          \lvert \partial_1 (t \epsilon S_v 
        + t^{\alpha +1} P_v) \rvert_{D_1} \leq t E_v, \
        \lvert \partial_2 (t \epsilon S_v + t^{\alpha +1} P_v) \rvert_{D_1} \leq t E_v.
    \end{equation*}
\end{corollary}

Similarly, let us use the fixed point theorem again to express $\Psi_v$ explicitly. Consider the second equation in \eqref{8}:
\begin{equation*}
    \theta = q  + \partial_1 \psi_v (I,q).
\end{equation*}
Let $X_2$ be the set of all real analytic functions $\theta_*(I,\theta)$ on $D_{r_v, s_v} $, where $\theta_*(I,\theta)$ is periodic $2 \pi$ with respect to angle variables, and$ \lvert \theta_*(I,\theta) \rvert _{r_v, s_v} \leq \sigma_v$.
Consider mapping:
\begin{equation*}
    f_2: \theta_* \to f_2(\theta_*) =  \partial_1 \psi_v(I, \theta - \theta_*),
\end{equation*}
We know that $f_2$ is well defined and real analytic on $X_2$, and $f_2 :X_2 \to X_2$. By Cauchy estimates, for $\forall \ \theta_* \in X_2, (I,\theta) \in D_{r_v, s_v}$, we have
\begin{equation*}
    \lvert f_2(\theta_*) \rvert \leq \lvert \partial_1 \psi_v(I, \theta - \theta_*) \rvert \leq \sigma_v.
\end{equation*}
For $\forall \ \theta_*^1, \theta_*^2 \in X_2$, $(I,\theta) \in D_{r_v, s_v}$, we have
\begin{align*}
    \lvert f_2(\theta_*^1) - f_2(\theta_*^2) \rvert  = &   \ \lvert \partial_1 \psi_v(I, \theta - \theta_*^1)
     -  \partial_1 \psi_v(I, \theta - \theta_*^2) \rvert \\
     \leq &  \ \frac{2 n c_1 \lvert \partial_1 \partial_2 \big( t \epsilon S_v(I, \theta - \bar{\theta}) + t^{\alpha +1} P_v(I, \theta - \bar{\theta})\big) \rvert}{t \gamma_v \sigma_v ^{\tau +n}}  \cdot \lvert p_*^1 - p_*^2 \rvert \\
     \leq & \ \frac{2 n c_1 E_v}{\gamma_v \sigma_v ^{\tau +n} \sigma_v } \cdot \lvert p_*^1 - p_*^2 \rvert \\
     \leq & \ L_2 \lvert p_*^1 - p_*^2 \rvert,
\end{align*}
if $\frac{2 n c_1 E_0}{\gamma_0 \sigma_0 ^{\tau +n+1} } <1 $, where $0<L_2<1 , \bar{\theta} \in X_2$.

So, $f_2: X_2 \to X_2$ is a compressed mapping, by the fixed point theorem, there is a unique fixed point $\theta_{**} \in X_2$, such that  $\Psi_v$ can be expressed explicitly as:
\begin{equation*}
    \begin{cases}
     p = I + \partial_2 \psi_v(I, \theta - \theta_{**}) ,\\
     q = \theta - \theta_{**},
    \end{cases}
\end{equation*}
so we get the following corollaries.

\begin{corollary} \label{corc2}
   $\partial_1 \psi_v(I,\theta - \theta_{**},\xi)$ and $\partial_2 \psi_v(I,\theta - \theta_{**},\xi)$ are well defined and real analytic on $D_2$, where $D_2= D_{r_v - 2 \rho_v, s_v - 2 \sigma_v} \times V_{\kappa^{'}_v} $, and it is periodic $2 \pi$ with respect to angle variables. Moreover, there are the following estimates on $D_2$:
   \begin{equation*}
          \lvert \partial_1 \psi_v \rvert_{D_2} \leq \frac{ 2 c_1 E_v}{\gamma_v \sigma_v ^{\tau +n}} , \ 
        \lvert \partial_2 \psi_v \rvert_{D_2} \leq \frac{2 c_1 E_v}{\gamma_v \sigma_v ^{\tau +n}}.
    \end{equation*}
\end{corollary}

\begin{corollary} \label{corc3}
    $\Psi_v$ and $\Psi_v^{-1}$ are well defined and real analytic on $D_3$, where $D_3= D_{r_v - 3 \rho_v, s_v - 3 \sigma_v} \times V_{\kappa^{'}_v} $, and they are both from $D_3$ to $ D_2$.
   Meanwhile, there is the following estimate on $D_3$:
    \begin{equation*}
       \lvert \Psi_v - id \rvert_{D_3} \leq \frac{2 c_1 E_v}{\gamma_v \sigma_v ^{\tau +n}} ,
    \end{equation*}
    so $\Psi_v^{-1} \circ R_v \circ C_v \circ \Psi_v$ is  well defined on $D_{r_{v+1}, s_{v+1}} \times V_{\kappa^{'}_v}  $, and from $D_{r_{v+1}, s_{v+1}}\times V_{\kappa^{'}_v}  $ to $D_1$.
\end{corollary}

Let $\Phi_{v+1} = \Phi_v \circ \Psi_v$, as we know
\begin{align*}
    D \Phi_v & = D (\Psi_0 \circ \Psi_1 \circ \dots \circ \Psi_{v-1}) \\
    & \leq \lvert D \Psi_0 \rvert  \lvert D \Psi_1 \rvert \dots \lvert D \Psi_{v-1} \rvert \\ 
    & \leq \prod \limits_v (1 + \frac{ 2 c_1 E_v}{\gamma_v \sigma_v ^{\tau +n}}) \\
    & \leq 2,
\end{align*}
if $\frac{2 c_1 E_0}{\gamma_0 \sigma_0 ^{\tau +n}} \leq 1$.
Then, on $D_{r_v-3\rho_v, s_v-3\sigma_v} $, we have
\begin{align*}
    \lvert \Phi_{v+1} - \Phi_v \rvert & =  \lvert \Phi_v \circ \Psi_v - \Phi_v \rvert \\
    & \leq \lvert  D \Phi_v \rvert \cdot  \lvert  \Psi_v - I \rvert \\
    & \leq \frac{4 c_1 E_v}{\gamma_v \sigma_v ^{\tau +n}} ,
\end{align*}
which means that $(\romannumeral1)$ in Lemma \ref{l1} (Iteration Lemma) holds. Besides,
\begin{align*}
    \lvert \omega_{v+1} - \omega_v \rvert & =  \big \lvert \big [ \partial_1  ( \epsilon S_v  + t^{\alpha} P_v )^* \big ] \big \rvert_{r_{v+1},s_{v+1}}  \leq 2 E_v,
\end{align*}
which means that $(\romannumeral3)$ in Lemma \ref{l1} (Iteration Lemma) holds.

Next, we proof $(\romannumeral2)$ in Lemma \ref{l1} (Iteration Lemma). Let $C_{v+1} = R_{v+1}^{-1} \circ \Psi_v^{-1} \circ R_v \circ C_v \circ \Psi_v$, $C_{v+1}$ is well defined and real analytic on $D_{r_{v+1},s_{v+1}} $, and then $C_{v+1}: (I, \theta) \to (\hat{I}, \hat{\theta})$  can be expressed as follows:
\begin{equation*}
    \begin{cases}
     \hat{I} = I + \partial_2 \psi_v (I,q) - \partial_2 \psi_v (\hat{I},\hat{q}) - \partial_2 \big(t \epsilon S_v
    + t^{\alpha +1} P_v\big)(\hat{p}, q), \\
     \hat{\theta} = \theta - \partial_1 \psi_v (I,q) + \partial_1 \psi_v (\hat{I},\hat{q}) + \partial_1 \big(t \epsilon S_v 
    + t^{\alpha +1} P_v\big)(\hat{p}, q)  - t \omega_{v+1} + t \omega_v .
    \end{cases}
\end{equation*}
Then, on $D_{r_{v+1},s_{v+1}} $, we have
\begin{align*}
    \lvert \hat{I} -  I \rvert = \ & \lvert  \partial_2 \psi_v (\hat{I},\hat{q}) - \partial_2 \psi_v (I,q)  + \partial_2 \big(t \epsilon S_v 
    + t^{\alpha +1} P_v \big)(\hat{p}, q) \rvert \\
     \leq \ & I_1 + I_2 + \dots + I_5 ,\\
    \lvert  \hat{\theta} -  \theta \rvert = \ & \lvert  \partial_1 \psi_v (\hat{I},\hat{q}) - \partial_1 \psi_v (I,q) + \partial_1 \big(t \epsilon S_v 
    + t^{\alpha +1} P_v \big)(\hat{p}, q) - t \omega_{v+1}+ t \omega_v \rvert \\
     \leq \ & J_1 + J_2 + \dots + J_5,
\end{align*}
where
\begin{align*}
    I_1 &= \lvert \partial_2 \psi_v (\hat{I},\hat{q}) - \partial_2 \psi_v (I,\hat{q}) \rvert ,\\
    I_2 &= \lvert \partial_2 \psi_v (I,\hat{q}) - \partial_2 \psi_v (I,q + t \omega_v) \rvert ,\\
    I_3 &= \lvert \partial_2 \psi_v (I,q + t \omega_v) - \partial_2 \psi_v (I,q) + \partial_2 \widetilde{\big(t \epsilon S_v(I,q)  + t^{\alpha +1} P_v(I,q) \big)^*} \rvert ,\\
    I_4 &= \lvert \partial_2 \widetilde{\big(t \epsilon S_v(\hat{p},q)  + t^{\alpha +1} P_v(\hat{p},q) \big)^*} - \partial_2 \widetilde{\big(t \epsilon S_v(I,q)  + t^{\alpha +1} P_v(I,q) \big)^*} \rvert, \\
    I_5 &= \lvert  \partial_2 \widetilde{\big(t \epsilon S_v(\hat{p},q)  + t^{\alpha +1} P_v(\hat{p},q) \big) }-  \partial_2 \widetilde{\big(t \epsilon S_v(\hat{p},q)  + t^{\alpha +1} P_v(\hat{p},q) \big)^*} \rvert, \\
    J_1 &= \lvert \partial_1 \psi_v (\hat{I},\hat{q}) - \partial_1 \psi_v (I,\hat{q}) \rvert, \\
    J_2 &= \lvert \partial_1 \psi_v (I,\hat{q}) - \partial_1 \psi_v (I,q + t \omega_v) \rvert, \\
    J_3 &= \lvert \partial_1 \psi_v (I,q + t \omega_v) - \partial_1 \psi_v (I,q) + \partial_1 \widetilde{\big(t \epsilon S_v(I,q)  + t^{\alpha +1} P_v(I,q) \big)^*} \rvert, \\
    J_4 &= \lvert \partial_1 \widetilde{\big(t \epsilon S_v(\hat{p},q)  + t^{\alpha +1} P_v(\hat{p},q) \big)^*} - \partial_1 \widetilde{\big(t \epsilon S_v(I,q)  + t^{\alpha +1} P_v(I,q) \big)^*} \rvert ,\\
    J_5 &= \lvert  \partial_1 \widetilde{ \big(t \epsilon S_v(\hat{p},q)  + t^{\alpha +1} P_v(\hat{p},q) \big)} -  \partial_1 \widetilde{\big(t \epsilon S_v(\hat{p},q)  + t^{\alpha +1} P_v(\hat{p},q) \big)^*} \rvert .
\end{align*}

\begin{lemma} \label{l6}
   (Lemma A.2 in \cite{r4}) If $f \in A^s$, then
   \begin{equation*}
       \lvert f - T_K f \rvert_{s-\sigma} \leq c K^n e^{-K \sigma} \lvert f \rvert_s ,   \   0 \leq \sigma \leq s,
   \end{equation*}
   where the constant $c$ only depends on $n$, $A^s$ is in Lemma \ref{l5}.
\end{lemma}

Next, by Lemma \ref{l5}, Lemma \ref{l6} and the Cauchy estimates, we estimate $I_1, I_2, \dots , I_5$:

\begin{align*}
    I_1 &= \lvert \partial_2 \psi_v (\hat{I},\hat{q}) - \partial_2 \psi_v (I,\hat{q}) \rvert_{r_v - 4 \rho_v, s_v - 4 \sigma_v} \\
    & \leq \frac{n \lvert \partial_2 \psi_v \rvert_{r_v - 3 \rho_v, s_v - 4 \sigma_v} }{\rho_v} \lvert \hat{I} - I \rvert_{r_v - 4 \rho_v, s_v - 4 \sigma_v} \\
    & \leq \frac{ 2 n c_1 E_v}{\gamma_v \sigma_v ^{\tau +n} \rho_v} \lvert \hat{I} - I \rvert_{r_v - 4 \rho_v, s_v - 4 \sigma_v}.
\end{align*}
\begin{align*}
    I_2 = & \ \lvert \partial_2 \psi_v (I,\hat{q}) - \partial_2 \psi_v (I,q + t \omega_v) \rvert_{r_v - 4 \rho_v, s_v - 4 \sigma_v} \\
     \leq & \ \frac{n \lvert \partial_2 \psi_v \rvert_{r_v - 4 \rho_v, s_v - 3 \sigma_v} }{\sigma_v} \lvert \hat{q} - q - t \omega_v \rvert_{r_v - 4 \rho_v, s_v - 4 \sigma_v} \\
     \leq & \ \frac{ 2 n c_1 E_v}{\gamma_v \sigma_v ^{\tau +n} \sigma_v} \cdot \lvert \hat{q} - q - t \omega_v \rvert_{r_v - 4 \rho_v, s_v - 4 \sigma_v} \\
     \leq & \ \frac{ 2 n c_1 E_v}{\gamma_v \sigma_v ^{\tau +n} \sigma_v} \cdot \lvert  \partial_1 \big(t \epsilon S_v(\hat{p}, q) + t^{\alpha +1} P_v(\hat{p}, q)\big) \rvert_{r_v - 4 \rho_v, s_v - 4 \sigma_v} \\
     \leq & \ \frac{ 2 n c_1 t  E_v^2}{\gamma_v \sigma_v ^{\tau +n+1} }.
\\
    I_3 = & \  \lvert \partial_2 \psi_v (I,q + t \omega_v) - \partial_2 \psi_v (I,q) + \partial_2 \widetilde{\big(t \epsilon S_v  + t^{\alpha +1} P_v \big)^*} (I,q) \rvert_{r_v - 4 \rho_v, s_v - 4 \sigma_v} \\
    \leq & \  \lvert \partial_2 \widetilde{\big(t \epsilon S_v  + t^{\alpha +1} P_v \big)^*}(I,q) - T_{K_v} \partial_2 \widetilde{\big(t \epsilon S_v  + t^{\alpha +1} P_v \big)^*} (I,q) \rvert_{r_v - 4 \rho_v, s_v - 4 \sigma_v}  \\ 
    \leq & \ c_2 K_v^n e^{-K_v \sigma_v} \lvert   \partial_2 \widetilde{\big(t \epsilon S_v(I,q)  + t^{\alpha +1} P_v(I,q) \big)^*} \rvert_{r_v - 4 \rho_v,s_v - 3 \sigma_v} \\
    \leq & \ 4 c_2 K_v^n e^{-K_v  \sigma_v} t E_v,
\end{align*}
just like Lemma \ref{l6}, here $c_2$ is a constant only depends on $n$.
\begin{align*}
    I_4 = & \ \lvert \partial_2 \widetilde{\big(t \epsilon S_v + t^{\alpha +1} P_v \big)^*}(\hat{p},q)  - \partial_2 \widetilde{\big(t \epsilon S_v + t^{\alpha +1} P_v \big)^*}(I,q) \rvert_{r_v - 4 \rho_v, s_v - 4 \sigma_v} \\
    \leq & \ \frac{n \lvert \partial_2 \widetilde{(t \epsilon S_v  + t^{\alpha +1} P_v )^*} \rvert_{r_v - 3 \rho_v, s_v - 4 \sigma_v}}{\rho_v} \lvert \hat{p} - I \rvert_{r_v - 4 \rho_v, s_v - 4 \sigma_v} \\
    \leq & \ \frac{4 n t E_v}{\rho_v} \{ \lvert \big(  \partial_2 \big(t \epsilon S_v(\hat{p}, q)  + t^{\alpha +1} P_v(\hat{p}, q)\big) \rvert + \lvert \partial_2 \psi_v (I,q) \rvert_{r_v - 4 \rho_v, s_v - 4 \sigma_v} \} \\
    \leq & \ \frac{4 n t E_v}{\rho_v} (t E_v + \frac{ 2 c_1   E_v}{\gamma_v \sigma_v ^{\tau +n} }) \\
    \leq & \ \frac{8 n c_1 t E_v^2}{\gamma_v \sigma_v ^{\tau +n} \rho_v}.
\\
    I_5 = & \ \lvert  \partial_2 \widetilde{\big(t \epsilon S_v  + t^{\alpha +1} P_v \big)}(\hat{p},q) 
     -  \partial_2 \widetilde{\big(t \epsilon S_v  + t^{\alpha +1} P_v \big)^*(\hat{p},q)} \rvert_{r_v - 4 \rho_v, s_v - 4 \sigma_v} \\
    \leq & \  \frac{\eta_v^2}{1-\eta_v} \lvert  \partial_2 \widetilde{\big(t \epsilon S_v(\hat{p},q)  + t^{\alpha +1} P_v(\hat{p},q) \big)} \rvert_{r_v , s_v} \\
    \leq & \ 4 \eta_v^2 t E_v, 
\end{align*}
if $\eta_v < \frac{1}{2}$.

In conclusion, if $\frac{ 2 n c_1 E_0}{\gamma_0 \sigma_0 ^{\tau +n} \rho_0} < 1 $, then we have
\begin{align*}
   & \  \lvert \hat{I} - I \rvert_{r_{v+1},s_{v+1}} \\
   \leq & \  \big( \  \frac{ 2 n c_1 }{\gamma_v \sigma_v ^{\tau +n+1} } + \frac{8 n c_1 }{\gamma_v \sigma_v ^{\tau +n} \rho_v} \ \big) t E_v^2 
      + ( \ 4 c_2 K_v^n e^{-K_v  \sigma_v} + 4 \eta_v^2 \  ) t E_v.
\end{align*}
Similarly, we can estimate $J_1, J_2, \dots , J_5$ the same way, then we have
\begin{align*}
    & \ \lvert \hat{\theta} - \theta \rvert_{r_{v+1},s_{v+1} } \\
    \leq & \  \big( \  \frac{ 2 n c_1 }{\gamma_v \sigma_v ^{\tau +n+1} } + \frac{8 n c_1 }{\gamma_v \sigma_v ^{\tau +n} \rho_v} \ \big) t E_v^2  + ( \ 4 c_2 K_v^n e^{-K_v  \sigma_v} + 4 \eta_v^2 \  ) t E_v.
\end{align*}

For $F_v = \frac{E_v}{\gamma_v \sigma_v^{ \tau +  n +1} \rho_v}$, if  $\epsilon$ and $t$ are small enough, we have $F_0$ is small enough such that $2n c_1 F_v \leq L_0, (0 < L_0 < 1)$, then
\begin{align*}
       & \ \frac{\big( \  \frac{ 2 n c_1 }{\gamma_v \sigma_v ^{\tau +n+1} } + \frac{8 n c_1 }{\gamma_v \sigma_v ^{\tau +n} \rho_v} \ \big) t E_v^2 + ( \ 4 c_2 K_v^n e^{-K_v  \sigma_v} + 4 \eta_v^2 \  ) t E_v
    }{\gamma_{v+1} \sigma_{v+1}^{ \tau +  n +1} \rho_{v+1}} \\
    \leq & \ \frac{c_3 \big \{\big( \  \frac{ 2 n c_1 }{\gamma_v \sigma_v ^{\tau +n+1} } + \frac{8 n c_1 }{\gamma_v \sigma_v ^{\tau +n} \rho_v} \ \big) t E_v^2 + ( \ 4 c_2 K_v^n e^{-K_v  \sigma_v} + 4 \eta_v^2 \  ) t E_v
    \big \}}{\gamma_v \sigma_v^{ \tau +  n +1} \eta_v \rho_v}  \\
    \leq & \ (10 n c_1 c_3) \cdot \frac{ t F_v^2}{\eta_v} + 
    ( 4 c_2 c_3 K_v^n e^{-K_v  \sigma_v} + 4 c_3 \eta_v^2 ) \cdot \frac{ t F_v}{\eta_v},
\end{align*}
where $c_3 = 2^{\tau +\bar{n} +n+2}$.
We let $\eta_v = \frac{1}{12 c_3} F_v^{\frac{1}{v}}$ for $v \geq 3$. If $K_0$ is big enough such that $144 c_2 c_3^2 K_v^n e^{-K_v  \sigma_v} \leq F_v^{\frac{2}{v}}$, and if $F_0$ is small enough such that $360 n c_1 c_3^2 \leq F_v ^{- (1-\frac{2}{v})}$ for $v \geq 3$, then we have
\begin{equation*}
     (10 n c_1 c_3) \cdot \frac{ t F_v^2}{\eta_v} +
    ( 4 c_2 c_3 K_v^n e^{-K_v  \sigma_v} + 4 c_3 \eta_v^2 ) \cdot \frac{ t F_v}{\eta_v}
    \leq t F_v^{\frac{v+1}{v}}, 
    \ v\geq 3  .
\end{equation*}
So we let $F_{v+1} =  F_v^{\frac{v+1}{v}}$, and by $ \lvert \partial_1 (t \epsilon S_v 
+ t^{\alpha +1} P_v) \rvert_{r_v,s_v} \leq t E_v$, $ \lvert \partial_2 (t \epsilon S_v + t^{\alpha +1} P_v) \rvert_{r_v,s_v} \leq t E_v$,
 we get that $C_{v+1}: (I, \theta) \to (\hat{I}, \hat{\theta})$
\begin{align*}
    & \lvert \hat{I} - I \rvert_{r_{v+1},s_{v+1}}, \ \lvert \hat{\theta} - \theta \rvert_{r_{v+1},s_{v+1}}\\
    \leq & \ \big( \  \frac{ 2 n c_1 }{\gamma_v \sigma_v ^{\tau +n+1} } + \frac{8 n c_1 }{\gamma_v \sigma_v ^{\tau +n} \rho_v} \ \big) t E_v^2 + ( \ 4 c_2 K_v^n e^{-K_v  \sigma_v} + 4 \eta_v^2 \  ) t E_v
 \end{align*}
 \begin{align*}
    \leq & \ \frac{ \big( \  \frac{ 2 n c_1 }{\gamma_v \sigma_v ^{\tau +n+1} } + \frac{8 n c_1 }{\gamma_v \sigma_v ^{\tau +n} \rho_v} \ \big) t E_v^2 + ( \ 4 c_2 K_v^n e^{-K_v  \sigma_v} + 4 \eta_v^2 \  ) t E_v }{\gamma_{v+1} \sigma_{v+1}^{ \tau +  n +1} \rho_{v+1} } \\
    & \ \cdot \gamma_{v+1} \sigma_{v+1}^{ \tau +  n +1} \rho_{v+1}  \\
    \leq & \ t F_v^{\frac{v+1}{v}} \cdot \gamma_{v+1} \sigma_{v+1}^{ \tau +  n +1} \rho_{v+1} \\
    \leq & \  t F_{v+1} \gamma_{v+1} \sigma_{v+1}^{ \tau +  n +1} \rho_{v+1} \\
    \leq & \ t E_{v+1}.
\end{align*}
Thus, $\lvert C_{v+1} - id \rvert _{r_{v+1}, s_{v+1}} \leq t E_{v+1}$.
This proofs $(\romannumeral2)$ in Lemma \ref{l1} (Iteration Lemma).

So, we complete the proof of the whole  Iteration Lemma.
$\square$

\subsection{Proof of convergence}

Note: As $F_{v+1} =  F_v^{\frac{v+1}{v}}$ for $v \geq 3$, we have $F_v = (((F_3^\frac{4}{3})^\frac{5}{4})^\frac{6}{5} \dots )^\frac{v}{v-1} = F_3^\frac{v}{3}$, then $\eta_v = \frac{1}{12 c_3} F_v^{\frac{1}{v}} = \frac{1}{12 c_3} F_3 ^\frac{1}{3}$, which is a constant.

Let $c_4 = \frac{1}{12 c_3} F_3 ^\frac{1}{3} = \frac{1}{3} \cdot 2^{-(\tau + \bar{n} +n +4)} F_3 ^\frac{1}{3}$, which is a constant, then by $r_v=\eta_{v-1} r_{v-1}$, we have $r_v=c_4 r_{v-1}$. In combination with the conditions mentioned earlier, we can get the final complete parameter settings:
\begin{equation*}
    \begin{array}{ccc}
    \tau_v = \frac{1^{-2}+\dots +v^{-2}}{2\sum \limits_{v=1}^\infty v^{-2}},
    & s_v = \frac{1}{4}(1-\tau_v) s_0 , \
    & \sigma_v = \frac{1}{4} (s_v - s_{v+1}) ,\\
    K_{v+1} = 4 K_v ,
    & r_{v+1}= c_4 r_v  , \
    & \rho_v = \frac{1}{4} (r_v - r_{v+1}) ,\\
     \gamma_v = \frac{\gamma^{\bar{n}+1}}{2^{(\bar{n}+1)v}},
    & F_v = \frac{E_v}{\gamma_v \sigma_v^{\tau+n+1}\rho_v}, 
    & F_{v+1} = F_v^{\frac{v+1}{v}},
\end{array}
\end{equation*}
where $c_4 = \frac{1}{3} \cdot 2^{-(\tau + \bar{n} +n +4)} F_3 ^\frac{1}{3}$.

Besides, as $\epsilon$ and $t$ are small enough, we can get that $F_0$(as well as $F_3$) is small enough and $K_0$ (as well as $K_3$) is big enough, then we can set appropriate initial values, so that the initial value conditions mentioned above are satisfied.

Going back to the Lemma \ref{l1} (Iteration Lemma), we have
\begin{align*}
     \lvert \Phi_v - id \rvert 
    & \leq \sum \limits_{v=0}^{\infty}  \lvert \Phi_{v+1} - \Phi_v \rvert \leq \sum \limits_{v=0}^{2} \lvert \Phi_{v+1} - \Phi_v \rvert +
    \sum \limits_{v=3}^{\infty} \frac{4 c_1  E_v}{\gamma_v \sigma_v^{\tau+n}} \\
    & \leq c_1^{'}( M_1 \epsilon +  2 n M_2 l_* t^{\alpha})   + 
    \sum \limits_{v=3}^{\infty} 4 c_1  F_v  \\ 
    & \leq c_1^{'}( M_1 \epsilon +  2 n M_2 l_* t^{\alpha})  + 4 c_1 \sum \limits_{v=3}^{\infty} F_3^{\frac{v}{3}} \\
    & \leq c_1^{'}( M_1 \epsilon +  2 n M_2 l_* t^{\alpha})   + 8 c_1 F_0 \\
    & \leq  c_2^{'}( M_1 \epsilon +  2 n M_2 l_* t^{\alpha}) ,
\\
     \lvert \omega_v - \omega \rvert & \leq \sum \limits_{v=0}^{\infty} \lvert \omega_{v+1} - \omega_v \rvert
    \leq \sum \limits_{v=0}^{2} \lvert \omega_{v+1} - \omega_v \rvert + \sum \limits_{v=3}^{\infty} 2 E_v \\
    & \leq  c_3^{'}( M_1 \epsilon +  2 n M_2 l_* t^{\alpha})    + 
    \sum \limits_{v=3}^{\infty} 2 F_v \\
    & \leq  c_3^{'}( M_1 \epsilon +  2 n M_2 l_* t^{\alpha})   + 2 \sum \limits_{v=3}^{\infty} F_3^{\frac{v}{3}} \\
    & \leq c_3^{'}( M_1 \epsilon +  2 n M_2 l_* t^{\alpha})   + 4  F_0  \\
    & \leq  c_4^{'}( M_1 \epsilon +  2 n M_2 l_* t^{\alpha}) , \\
    \lvert C_v - id \rvert & \leq t E_v \leq t F_v = t (F_3^{\frac{1}{3}})^v
    \leq t \big(c_5^{'}( M_1 \epsilon +  2 n M_2 l_* t^{\alpha})  \big)^v,
\end{align*}
for $v= 0,1,2,3, \dots$, where $c_2^{'}, c_4^{'}$ and $c_5^{'}$ depend only on $\tau, n, \bar{n},  \rho_0, \sigma_0, \gamma_0 $ and $K_0$.
So, as $\epsilon$ and $t$ are small enough, all of them are small enough.

Let $\Phi_{\epsilon,t} = \lim \limits_{v \to \infty} \Phi_v$,
$\omega_{\epsilon,t} = \lim \limits_{v \to \infty} \omega_v$, $R_{\epsilon,t}(p,q) = (p, q + t \omega_{\epsilon,t} )$ and
 $V_{\epsilon,t} = \bigcap \limits_{v} V_{\gamma,t,v}$, then
$\Phi _{\epsilon,t}^{-1} \circ G_{H^ \epsilon}^t \circ \Phi _{\epsilon,t} = R_{\epsilon,t}$, where $R_{\epsilon,t}$ is a rotation on $V_{\epsilon,t} \times \mathbb{T}^n$ with frequency $t \omega_{\epsilon,t} $, so
\begin{align*}
    & \lvert \Phi_{\epsilon,t} - id \rvert \leq  c_2^{'}( M_1 \epsilon +  2 n M_2 l_* t^{\alpha}), \\
    & \lvert \omega_{\epsilon,t} - \omega \rvert \leq  c_4^{'}( M_1 \epsilon +  2 n M_2 l_* t^{\alpha}).
\end{align*}

Let the time step be $t_1$ and $t_2$ respectively, and the corresponding symbols of the system be added $t_1$ and $t_2$ as superscripts, similarly constructing $\psi^{t_1}_{v}$ and $\psi^{t_2}_{v}$, then for $\xi \in V_{\epsilon,t_1} \cap V_{\epsilon,t_2}$ we have
\begin{align*}
    & \lvert  \partial_1 \psi^{t_1}_{v}  \rvert _{s_v - \sigma_v}  \leq \frac{2 c_1}{t_1 \gamma_v \sigma_v ^{\tau +n}} \lvert \partial_1  (t_1 \epsilon H_v  + t_1^{\alpha +1} P_v ) \rvert_{s_v} , \\
    & \lvert  \partial_2 \psi^{t_1}_{v}  \rvert _{s_v - \sigma_v}  \leq \frac{2 c_1}{t_1 \gamma_v \sigma_v ^{\tau +n}} \lvert \partial_2  (t_1 \epsilon H_v  + t_1^{\alpha +1} P_v ) \rvert_{s_v} , \\
    & \lvert  \partial_1 \psi^{t_2}_{v}  \rvert _{s_v - \sigma_v}  \leq \frac{2 c_1}{t_2 \gamma_v \sigma_v ^{\tau +n}} \lvert \partial_1  (t_2 \epsilon H_v  + t_2^{\alpha +1} P_v ) \rvert_{s_v} , \\
     & \lvert  \partial_2 \psi^{t_2}_{v}  \rvert _{s_v - \sigma_v}  \leq \frac{2 c_1}{ t_2 \gamma_v \sigma_v ^{\tau +n}} \lvert \partial_2  (t_2 \epsilon H_v  + t_2^{\alpha +1} P_v ) \rvert_{s_v} .
\end{align*}
Comparing these two systems, combined with corollary \ref{corc2} and corollary \ref{corc3}, we have 
\begin{equation*}
          \lvert  \Psi_v^{t_1} -  \Psi_v^{t_2} \rvert_{D_3} \leq \frac{ 2 c_1 \cdot 2 n M_2 l_* (t_1^{\alpha}-t_2^{\alpha})}{\gamma_v \sigma_v ^{\tau +n}} .
    \end{equation*}
So, 
\begin{equation*}
    \lvert \Phi_v^{t_1} - \Phi_v^{t_2} \rvert \leq  2 c_2^{'}  n M_2 l_* (t_1^{\alpha}-t_2^{\alpha}) , \ \lvert \omega_v^{t_1} - \omega_v^{t_2} \rvert \leq  2 c_4^{'} n M_2 l_* (t_1^{\alpha}-t_2^{\alpha}),
\end{equation*}
thus 
    \begin{equation*}
    \lvert \Phi_{\epsilon,t_1} - \Phi_{\epsilon,t_2} \rvert \leq  2 c_2^{'}  n M_2 l_* (t_1^{\alpha}-t_2^{\alpha})  , \ \lvert \omega_{\epsilon,t_1} - \omega_{\epsilon,t_2} \rvert \leq  2 c_4^{'} n M_2 l_* (t_1^{\alpha}-t_2^{\alpha}),
\end{equation*}
where $M_2$ and $l_*$ are constants, and $c_2^{'}$ and $c_4^{'}$ are constants depending only on $\tau, n, \bar{n},  \rho_0, \sigma_0, \gamma_0 $ and $K_0$.

As for $(\romannumeral2)$ in theorem \ref{t1} , we place the measure estimations in the next subsection.

\subsection{Measure Estimation}
\label{s34}

Let $V_{\epsilon,t} = \bigcap \limits_{v} V_{\gamma,t,v}$, where
\begin{equation*}
    V_{\gamma,t,v} = \big \{\xi \in V : \lvert e^{i \langle k, t \omega_v(\xi) \rangle} - 1 \rvert \geq \frac{t \gamma_v}{{\lvert k \rvert}^\tau} , \forall \  k \in \mathbb{Z}^n \setminus \{0 \} \big \}.
\end{equation*}
We focus on
\begin{equation*}
     \big \{\xi \in V : \lvert  \langle k, t \omega_v(\xi) \rangle -2 \pi l  \rvert \geq \frac{t \gamma_v}{{\lvert k \rvert}^\tau} , \forall \  k \in \mathbb{Z}^n \setminus \{0 \}, \forall  \ l \in \mathbb{Z} \big \}.
\end{equation*}
It is easy to know that the latter is contained in the former.
 Before estimating its  measure, we first introduce the following lemma.

\begin{lemma} \label{l8}
   (Lemma 4.9 in \cite{r6} ) Let $ K \subseteq \mathbb{R}^n $ be compact with positive diameter $d := \sup_{x,y \in K} | x-y|_2 >0$,  define $\mathcal{B}:= ({K +\theta }) \cap \mathbb{R}^n \subseteq \mathbb{R }^ n $ for some $\theta > 0$,
   and $ g \in C^{u_0+1}(\mathcal{B},\mathbb{R})$ be a function with
\begin{equation} \label{31}
    \min_{y\in K} \max_{0\leq v\leq {u_0}} \lvert D^v g(y) \rvert \geq \beta,
\end{equation}
   for some ${u_0} \in \mathbb{N} $ and $\beta>0$. Then for any $\tilde g\in C^{u_0}(\mathcal{B},\mathbb{R})$ satisfying $\lvert \tilde g - g  \rvert^{u_0}_{\mathcal{B}} := \max \limits_{0\leq v\leq {u_0}} |D^v (\tilde g - g)|_{\mathcal{B}} \leq \frac{1}{2} \beta$, we have the estimate
\begin{equation} \label{32}
    \big \lvert  \{ y\in K : \lvert \tilde g(y) \rvert \leq \epsilon  \} \big \rvert \leq B d^{n-1}(n^{-\frac{1}{2}}+2d+\theta^{-1}d)(\frac{\epsilon}{\beta})^\frac{1}{u_0}\frac{1}{\beta}\max_{0 < v\leq {u_0}} \lvert D^v g \rvert_{\mathcal{B}},
\end{equation}
    whenever $0<\epsilon \leq  \frac{\beta}{2u_0 +2}$. Here, $B =3(2\pi e)^\frac{n}{2} (u_0+1)^{u_0+2} [(u_0 +1)!]^{-1}$.
\end{lemma}

\begin{remark}
   This Lemma is just the Lemma 4.9 in \cite{r6}, and it comes from Theorem 17.1 in \cite{r5} given by R\"{u}ssmann. The proof is similar to that.
\end{remark}

Firstly, combined with Remark \ref{re1}, we know that $t \omega$ satisfies the R\"{u}ssmann's non-degeneracy condition, and by (\ref{30}), there exist $\bar{n} =\bar{n}(\omega, V) \in \mathbb{N}$ and $\beta = \beta(\omega, V) >0$ such that $\min \limits_{\xi \in V} \max \limits_{0\leq v\leq \bar{n}} |D^v \langle k, t \omega(\xi) \rangle |\geq t \beta \lvert k \rvert$ for $\forall \ k \in \mathbb{Z}^n \setminus \{0 \}$.

As for $ \lvert \langle k, t\omega \rangle -2 \pi l  \rvert $, let $\widetilde{V} = V \times (0,1)$, $\widetilde{\xi} = (\xi, \xi^{'})$, $\xi^{'} \in (0,1)$, $t \widetilde{ \omega} = (t \omega, -2 \pi)$, $\widetilde{k} = (k, l)$, then we have $ rank \big \{\partial^i_{\xi} \widetilde{\omega}(\widetilde{\xi}) :  |i| \leq \bar{n}+1 \big \} = n+1 ,  \   \widetilde{\xi} \in \widetilde{V}$, and
$\min \limits_{\widetilde{\xi} \in \widetilde{V}} \max \limits_{0\leq v\leq \bar{n}+1} |D^v \langle \widetilde{k}, t \widetilde{\omega}(\widetilde{\xi}) \rangle |\geq t \beta \lvert \widetilde{k} \rvert$.

Let $\mathcal{B} = (\widetilde{V}+ \kappa) \cap \mathbb{R}^{n+1}$, as
$\lvert \omega_v - \omega \rvert  \leq c_4^{'}( M_1 \epsilon +  2 n M_2 l_* t^s)  $,
then for $\epsilon$ and $t$ small enough, we have
$ \lvert \langle \widetilde{k}, t \widetilde{\omega}_v \rangle  - \langle \widetilde{k}, t \widetilde{\omega}  \rangle \rvert^{\bar{n}+1}_{\mathcal{B}}  \leq \frac{1}{2} t \beta \lvert \widetilde{k} \rvert
$. And for $\gamma_0$ small enough, we have $t\gamma_v \lvert k\rvert^{-\tau} \leq \frac{t \beta \lvert \widetilde{k} \rvert}{2 (\bar{n}+1) +2}$, then we can use Lemma \ref{l8}, and we have
\begin{align*}  \label{33}
    & \big \lvert  \{ \xi \in V : \lvert \langle k, t\omega_v \rangle - 2 \pi l   \rvert < t\gamma_v \lvert k\rvert^{-\tau}  \} \big \rvert =  \big \lvert  \{ \widetilde{\xi} \in \widetilde{V} : \lvert \langle \widetilde{k}, t \widetilde{\omega}_v \rangle \rvert < t\gamma_v \lvert k\rvert^{-\tau}  \} \big \rvert
    \\
    & \ \leq \ c_5 d^{n}((n+1)^{-\frac{1}{2}}+2d+ \kappa^{-1}d) (\frac{t \gamma_v}{t \beta})^\frac{1}{\bar n+1} \lvert k \rvert^{\frac{ -\tau } {\bar{n} +1} } \lvert \widetilde{k} \rvert^{\frac{ -1 } {\bar{n} +1} }  \frac{1}{t \beta}  \max_{0 < v\leq \bar{n}+2} \lvert D^v t \widetilde{\omega} \rvert_{\mathcal{B}},
\end{align*}
for $  k\in \mathbb{Z}^n \setminus \{0\} ,  l \in \mathbb{Z}$ , where $ c_{5} = 3(2\pi e)^\frac{n+1}{2} (\bar n +2)^{\bar n +3} [(\bar n  +2)!]^{-1} $, $d$ is the diameter of $V \times (0,1)$.
Note that for such $l$, we have $\lvert l \rvert \leq \lvert k \rvert  \hat{K} $, where $\hat{K} = \lvert \omega \rvert_{V} + c_4^{'}( M_1  +  2 n M_2 l_* )   +1$.
Define
 \begin{align*}
    R^{k,l}_v &= \big \{ \xi \in V :  \lvert \langle k, t\omega_v \rangle -2 \pi l   \rvert  < t\gamma_v \lvert k \rvert^{-\tau}  \big \}, \ \lvert l \rvert \leq \lvert k \rvert  \hat{K}  , \  k\in \mathbb{Z}^n \setminus \{0\} , \\ 
    R^k_v &= \big \{ \xi \in V :  \lvert \langle k, t\omega_v \rangle -2 \pi l  \rvert  < t\gamma_v \lvert k \rvert^{-\tau} , \exists \ l \in \mathbb{Z}   \big \}, \  k\in \mathbb{Z}^n \setminus \{0\} ,\\
    R_v &= \big \{ \xi \in V :  \lvert \langle k, t\omega_v \rangle  -2 \pi l  \rvert  < t\gamma_v \lvert k \rvert^{-\tau} ,\exists \ l \in \mathbb{Z}  ,  \exists  \ k\in \mathbb{Z}^n \setminus \{0\} \big \},
\end{align*}
i.e., $ R_v= \bigcup \limits_{k\in \mathbb{Z}^n \setminus \{0\}} R^k_v = \bigcup \limits_{k\in \mathbb{Z}^n \setminus \{0\} } \bigcup \limits_{\lvert l \rvert \leq \lvert k \rvert  \hat{K} }  R^{k,l}_v$.
Besides,  $\lvert \widetilde{k} \rvert^{\frac{ -1 } {\bar{n} +1} } \leq 1$ for all $ k\in \mathbb{Z}^n \setminus \{0\}$. Then, we have
\begin{align*}
    \lvert R^{k,l}_v \rvert & \leq c_5 d^{n}((n+1)^{-\frac{1}{2}}+2d+ \kappa^{-1}d) (\frac{t \gamma_v}{t \beta})^\frac{1}{\bar n+1}  \frac{1}{t \beta} \max_{0 < v\leq \bar{n}+2} \lvert D^v t \widetilde{\omega} \rvert_{\mathcal{B}} \lvert k \rvert^{\frac{ -\tau} {\bar n +1 }} \\
    & \leq c_5 d^{n}((n+1)^{-\frac{1}{2}}+2d+ \kappa^{-1}d) (\frac{ \gamma_v}{ \beta})^\frac{1}{\bar n+1}  \frac{1}{ \beta} \max_{0 < v\leq \bar{n}+2} \lvert D^v  \omega \rvert_{\mathcal{B}} \lvert k \rvert^{\frac{ -\tau} {\bar n +1 }} ,
\\
    \lvert R^k_v \rvert & = \lvert \bigcup \limits_{\lvert l \rvert \leq \lvert k \rvert  \hat{K} }  R^{k,l}_v \rvert  \leq 2 \lvert k \rvert \hat{K} \lvert R^{k,l}_v \rvert \\
    & \leq 2 \hat{K} c_5   d^{n}((n+1)^{-\frac{1}{2}}+2d+ \kappa^{-1}d) (\frac{ \gamma_v}{ \beta})^\frac{1}{\bar {n}+1}  \frac{1}{ \beta} \lvert  \omega \rvert^{\bar{ n}  +2}_{\mathcal{B}} \lvert k \rvert^{\frac{\bar{n} +1 -\tau} {\bar {n} +1 }} \\
    & \leq c_6   d^{n}((n+1)^{-\frac{1}{2}}+2d+ \kappa^{-1}d) (\frac{ \gamma_v}{ \beta})^\frac{1}{\bar {n}+1}  \frac{1}{ \beta} \lvert \omega \rvert^{\bar{ n}  +2}_{\mathcal{B}} \lvert k \rvert^{\frac{\bar{n} +1 -\tau} {\bar {n} +1 }} ,
\end{align*}
where $c_6 = 6 (2\pi e)^\frac{n+1}{2} (\bar n +2)^{\bar n +3} [(\bar n  +2)!]^{-1} \hat{K}$.

By $V_{\epsilon,t} = \bigcap \limits_{v=0}^{\infty} V_{\gamma,t,v}$,
we have the following inequality.
\begin{align*}
    & \lvert V \setminus V_{\epsilon,t} \rvert  \leq \sum _{v=0}^\infty \lvert R_v \rvert  \ \leq \sum _{v=0}^\infty  \sum _{k\in \mathbb{Z}^n \setminus \{0\} } \lvert R_v^k \rvert 
    \leq \sum _{v=0}^\infty \sum^\infty _{r=1} 2n(2r+1)^{n-1}\lvert R_v^r \rvert \\
   & \leq \sum _{v=0}^\infty \sum^\infty _{r=1} 2n(2r+1)^{n-1} c_6   d^{n}((n+1)^{-\frac{1}{2}}+2d+ \kappa^{-1}d) (\frac{ \gamma_v}{ \beta})^\frac{1}{\bar {n}+1}  \frac{1}{ \beta} \lvert  \omega \rvert^{\bar{ n}  +2}_{\mathcal{B}} r^{\frac{\bar{n} +1 -\tau} {\bar {n} +1 }} \\
   &\leq \sum ^\infty _{r=1} 2n(2r+1)^{n-1} r^{\frac{\bar{n} +1 -\tau} {\bar {n} +1 }}  c_6   d^{n}((n+1)^{-\frac{1}{2}}+2d+ \kappa^{-1}d) \\
     &  \ \cdot \sum _{v=0}^\infty \gamma_v^{\frac{1}{\bar {n}+1}}  \beta^{-\frac{\bar{n}+2}{\bar{n } +1}} \lvert  \omega \rvert^{\bar{ n}  +2}_{\mathcal{B}} .
\end{align*}
And for $\tau \geq (n+2)(\bar{n}+1) $,
$\bar{c} := \sum \limits_{r=1}^\infty  2n(2r+1)^{n-1} r^{\frac{\bar{n} +1 -\tau} {\bar {n} +1 }} $ is convergent.
Besides,  $\gamma_v =  \frac{\gamma ^{\bar{n}+1}}{2^{(\bar n +1)v}} $,
we have $ \sum \limits_{v=0}^\infty \gamma_v^{\frac{1}{\bar n+1}} = \sum \limits_{v=0}^\infty \frac{\gamma}{2^v} = 2 \gamma $, and $\lvert  \omega \rvert^{\bar{ n}  +2}_{\mathcal{B}}  < \infty$ by analysis, therefore,
\begin{align*}
    \lvert V \setminus V_{\epsilon,t} \rvert & \leq \bar{c} \cdot 2\gamma \cdot   c_6  d^{n}((n+1)^{-\frac{1}{2}}+2d+ \kappa^{-1}d)  \beta^{-\frac{\bar{n}+2}{\bar{n } +1}} \lvert  \omega \rvert^{\bar{ n}  +2}_{\mathcal{B}}  \\
    & \leq c_7   \gamma  d^{n}((n+1)^{-\frac{1}{2}}+2d+ \kappa^{-1}d),
\end{align*}
where $c_7 = 12 \bar{c} (2\pi e)^\frac{n+1}{2} (\bar n +2)^{\bar n +3} [(\bar n  +2)!]^{-1} \hat{K} \beta^{-\frac{\bar{n}+2}{\bar{n } +1}} \lvert  \omega \rvert^{\bar{ n}  +2}_{\mathcal{B}} $,
$\bar{c}$ depends only on $n$, $\bar{n}$, $\tau$.
So for $\gamma$ small enough ,  $V_{\epsilon,t}$ can be a set of positive measures. And for $\gamma \to 0$, we have $\lvert V \setminus V_{\epsilon,t} \rvert \to 0 $. Then
$(\romannumeral2)$ in theorem \ref{t1} proves.

Similarly, we know that $ V_{\epsilon,t_1} \cap V_{\epsilon,t_2} $ is also a set of positive measure if $\gamma$ is small enough. Since  $\lvert V \setminus (V_{\epsilon,t_1} \cap V_{\epsilon,t_2} ) \rvert \leq \lvert V \setminus V_{\epsilon,t_1}  \rvert + \lvert V \setminus  V_{\epsilon,t_2} \rvert  $, we also have $\lvert V \setminus (V_{\epsilon,t_1} \cap V_{\epsilon,t_2}) \rvert \to 0 $ as $\gamma \to 0$ .

The proof of Theorem \ref{t1} is complete. $\square$

\section{Kolmogorov's non-degeneracy condition}

\begin{remark} \label{rem2}
   Without the R\"{u}ssmann's non-degeneracy condition, Theorem \ref{t1} also holds for $\omega$ satisfying the Kolmogorov's non-degeneracy condition $\Theta_1 \lvert \xi_1 - \xi_2 \rvert \leq \lvert \omega(\xi_1) - \omega(\xi_2) \rvert \leq \Theta_2 \lvert \xi_1 - \xi_2 \rvert$, where $\Theta_1$ and $\Theta_2$ are constants.
   Moreover, when time step  $t$ and disturbance parameter  $\epsilon$ are small enough, 
both the generating function representation of symplectic algorithm \eqref{2} and 
the generating function representation of phase flow in 
nearly integrable Hamiltonian system \eqref{xy3}  have invariant tori on corresponding sets of large  measures, and there is only a slight deformation in the common area. The Hausdorff's distance between the invariant tori in the phase spaces of the two systems is about $o(t^{\alpha})$.
\end{remark}

We explain this remark in three steps.

Step1:

   To be specific, if $\tau > n+2$, $\gamma_v$ is changed to $\frac{\gamma}{2^{v}}$, $c_3$ and $c_4$ are changed to $2^{\tau +n +2}$ and $ \frac{1}{3} \cdot 2^{-(\tau +n +4)} F_3 ^\frac{1}{3}$ respectively, and the other parameters reduce the dependence on $\bar{n}$, then we can construct a similar iterative procedure.
That is to say, there is a non-empty Cantor set $V_{\epsilon,t} \subseteq V$ and a Whitney smooth symplectic mapping  
$\Phi _{\epsilon,t}  $ such that 
   $\Phi _{\epsilon,t}^{-1} \circ G_{H^ \epsilon}^t \circ \Phi _{\epsilon,t} = R_{\epsilon,t}$, and 
\begin{align*}
    & \lvert \Phi_{\epsilon,t} - id \rvert \leq  c_2^{'}( M_1 \epsilon +  2 n M_2 l_* t^{\alpha}), \\
    & \lvert \omega_{\epsilon,t} - \omega \rvert \leq  c_4^{'}( M_1 \epsilon +  2 n M_2 l_* t^{\alpha}),
\end{align*}
 where $c_2^{'}$ and $ c_4^{'}$ depend only on $\tau, n, \rho_0, \sigma_0, \gamma_0 $ and $K_0$.

   As for the corresponding measure estimation, let $\frac{t \gamma}{{\lvert k \rvert}^\tau} <1$, then we have
   \begin{equation*}
        |R^{k,l}_v| = |\big \{ \xi \in V :  \lvert \langle k, t\omega_v \rangle -2 \pi l   \rvert  < t\gamma_v \lvert k \rvert^{-\tau}  \big \}| \leq \frac{2 \gamma_v}{|k|^{\tau +1} \Theta_1}.
   \end{equation*}
   So,
   \begin{equation*}
        \lvert R^k_v \rvert  \leq 2 \lvert k \rvert \hat{K} \lvert R^{k,l}_v \rvert \leq \frac{4 \gamma_v \hat{K}}{|k|^{\tau} \Theta_1},
   \end{equation*}
   here $\hat{K}$ stays the same, then, we get
   \begin{align*}
     \lvert V \setminus V_{\epsilon,t} \rvert & \leq \sum _{v=0}^\infty \lvert R_v \rvert  \ \leq \sum _{v=0}^\infty  \sum _{k\in \mathbb{Z}^n \setminus \{0\} } \lvert R_v^k \rvert 
    \leq \sum _{v=0}^\infty \sum^\infty _{r=1} 2n(2r+1)^{n-1}\lvert R_v^r \rvert \\
   & \leq \sum _{v=0}^\infty \sum^\infty _{r=1} 2n(2r+1)^{n-1} \frac{4 \gamma_v \hat{K}}{r^{\tau} \Theta_1} \leq \sum^\infty _{r=1} 2n(2r+1)^{n-1} \frac{8 \gamma_0 \hat{K}}{r^{\tau} \Theta_1} \\
   & \leq \bar{c} \frac{8 \gamma_0 \hat{K}}{ \Theta_1}, 
\end{align*}
where $\bar{c} := \sum \limits_{r=1}^\infty  2n(2r+1)^{n-1} \frac{1}{r^{\tau}} $ is convergent.
That is to say, we also get that $V_{\epsilon,t}$ can be a set of positive measures for $\gamma$ small enough, and $\lvert V \setminus V_{\epsilon,t} \rvert \to 0 $ as $\gamma \to 0$. 

   And, more importantly, in the Kolmogorov's non-degeneracy condition, the frequency mapping is a local differential homeomorphism. For all $w_*=\omega(\xi), \xi \in V \setminus V_{\epsilon,t}$, there exist $x_v \in V$, such that $\omega_v(x_v) = w_*$, $v=0,1,2, \dots$. Let $v \to \infty$, then there exists $x_0 \in V$, such that $\omega_0(x_0) = w_* = \omega_{\infty}(x_{\infty})$, and we have 
\begin{equation*}
       \lvert x_{\infty} - x_0 \rvert \leq \frac{c_4^{'}( M_1 \epsilon +  2 n M_2 l_* t^{\alpha})}{\Theta_1}.
\end{equation*}

Step 2:

We can also do the same for the generating function representation of phase flow in 
nearly integrable Hamiltonian system \eqref{xy3}. 
For $\xi \in \widetilde{V}_{\epsilon,t}$, where 
$\widetilde{V}_{\epsilon,t} = \bigcap \limits_{v=0}^\infty \widetilde{V}_{\gamma,t,v} $, we have
\begin{align*}
     \lvert \widetilde{\Phi}_v - id \rvert & \leq  \sum \limits_{v=0}^{\infty}  \lvert \widetilde{\Phi}_{v+1} - \widetilde{\Phi}_v \rvert \leq \sum \limits_{v=0}^{2} \lvert \widetilde{\Phi}_{v+1} - \widetilde{\Phi}_v \rvert +
    \sum \limits_{v=3}^{\infty} \frac{4 c_1  \widetilde{E}_v}{\gamma_v \sigma_v^{\tau+n}} \\
    & \leq \widetilde{c}_1 \cdot M_1 \epsilon    + 8 c_1 \widetilde{F}_0
     \leq  \widetilde{c}_2\cdot M_1 \epsilon, \\
    \lvert \widetilde{\omega}_v - \omega  \rvert & \leq \sum \limits_{v=0}^{\infty} \lvert \widetilde{\omega}_{v+1} - \widetilde{\omega}_v \rvert
    \leq \sum \limits_{v=0}^{2} \lvert \widetilde{\omega}_{v+1} - \widetilde{\omega}_v \rvert + \sum \limits_{v=3}^{\infty} 2 \widetilde{E}_v \\
    & \leq \widetilde{c}_3 \cdot M_1 \epsilon  + 4  \widetilde{F}_0   \leq  \widetilde{c}_4 \cdot M_1 \epsilon,
\end{align*}
where $\widetilde{c}_2$ and $\widetilde{c}_4$ depend only on $\tau, n,  \rho_0, \sigma_0, \gamma_0 $ and $K_0$.
Let $\widetilde{\Phi}_{\epsilon,t} = \lim \limits_{v \to \infty} \widetilde{\Phi}_v$,
$\widetilde{\omega}_{\epsilon,t} = \lim \limits_{v \to \infty} \widetilde{\omega}_v$ , then we have
\begin{equation*}
    \lvert \widetilde{\Phi}_{\epsilon,t} - id \rvert \leq  \widetilde{c}_2 \cdot M_1 \epsilon,  \
     \lvert \widetilde{\omega}_{\epsilon,t} - \omega \rvert \leq  \widetilde{c}_4 \cdot M_1 \epsilon.
\end{equation*}

Step 3:

What is more, since the frequency mapping is a local differential homeomorphism, we can compare the generating function representation of symplectic algorithm \eqref{2} and the generating function representation of phase flow in 
nearly integrable Hamiltonian system \eqref{xy3}.

For the time step $t$ in the common area, we have
\begin{equation*} 
    \psi_v - \widetilde{\psi}_v 
    = - \sum \limits_{0< \lvert k \rvert \leq K_v} \frac{\big(  ( t \epsilon S_v - t \epsilon \widetilde{S}_v )^* \big )_k e^{i \langle k, \theta   \rangle} + \big(  ( t^{\alpha +1} P_v )^* \big )_k e^{i \langle k, \theta \rangle}} {e^{i \langle k, t \omega_v \rangle} -1} ,
\end{equation*}
and 
\begin{align} 
   & \lvert \partial_1 \psi_v - \partial_1 \widetilde{\psi}_v \rvert _{s_v - \sigma_v} \leq \frac{2 c_1}{ \gamma_v \sigma_v ^{\tau +n}} \lvert \partial_1  (\epsilon S_v - \epsilon \widetilde{S}_v + t^{\alpha } P_v ) \rvert_{s_v} , \notag \\
   & \lvert \partial_2 \psi_v - \partial_1 \widetilde{\psi}_v \rvert _{s_v - \sigma_v} \leq \frac{2 c_1}{ \gamma_v \sigma_v ^{\tau +n}} \lvert \partial_2  (\epsilon S_v - \epsilon \widetilde{S}_v + t^{\alpha } P_v ) \rvert_{s_v} . 
\end{align}
Let $F_v^* = \frac{E_v^*}{\gamma_v \sigma_v^{\tau+n+1}\rho_v} ,
     F_{v+1}^* = (F_v^*)^{\frac{v+1}{v}}$,
$E_0^* = E_0 - \widetilde{E}_0$, where $E_0^* = o(t^{\alpha})$, then  for the time step $t$ in the common area, combined with \eqref{xy15}, we have
\begin{align*}
    & \lvert (\Phi_v - \widetilde{\Phi}_v) - (\Phi_{v-1} -\widetilde{\Phi}_{v-1} ) \rvert \leq \frac{4 c_1  E_{v-1}^*}{\gamma_{v-1} \sigma_{v-1}^{\tau+n}} ,  \\ 
    & \lvert (\omega_v - \widetilde{\omega}_v) - (\omega_{v-1} - \widetilde{\omega}_{v-1}) \rvert \leq 2 E_v^*.
\end{align*}
Then, for the common area $\xi \in \bigcap \limits_{v=0}^\infty V_{\gamma,t,v} \cap \widetilde{V}_{\gamma,t,v}$, we have 
\begin{align*}
     \lvert \Phi_v - \widetilde{\Phi}_v \rvert & \leq  \sum \limits_{v=0}^{\infty}  \lvert (\Phi_{v+1} - \widetilde{\Phi}_{v+1}) - (\Phi_v - \widetilde{\Phi}_v) \rvert \\
     & \leq \sum \limits_{v=0}^{2} \lvert (\Phi_{v+1} - \widetilde{\Phi}_{v+1}) - (\Phi_v - \widetilde{\Phi}_v) \rvert +
    \sum \limits_{v=3}^{\infty} \frac{4 c_1  E_v^*}{\gamma_v \sigma_v^{\tau+n}} \\
    & \leq c_1^* \cdot  t^{\alpha}    + 8 c_1 F_0^*
     \leq  c_2^* \cdot  t^{\alpha}, \\
    \lvert \omega_v - \widetilde{\omega}_v  \rvert & \leq \sum \limits_{v=0}^{\infty} \lvert (\omega_{v+1} - \widetilde{\omega}_{v+1}) - (\omega_v - \widetilde{\omega}_v) \rvert \\
    & \leq \sum \limits_{v=0}^{2} \lvert (\omega_{v+1} - \widetilde{\omega}_{v+1}) - (\omega_v - \widetilde{\omega}_v) \rvert + \sum \limits_{v=3}^{\infty} 2 E_v^* \\
    & \leq c_3^* \cdot  t^{\alpha}  + 4  F_0^*   \leq  c_4^* \cdot  t^{\alpha},
\end{align*}
where $c_2^*$ and $c_4^*$ depend only on $\tau, n,   \rho_0, \sigma_0, \gamma_0 $ and $K_0$. And by $\Phi_{\epsilon,t} - \widetilde{\Phi}_{\epsilon,t} = \lim \limits_{v \to \infty}  \Phi_v - \widetilde{\Phi}_v  $, $\omega_{\epsilon,t} - \widetilde{\omega}_{\epsilon,t} = \lim \limits_{v \to \infty} \omega_v - \widetilde{\omega}_v$, we have
\begin{equation*}
     \lvert \Phi_{\epsilon,t} - \widetilde{\Phi}_{\epsilon,t} \rvert  \leq c_2^* t^{\alpha}, \ \lvert \omega_{\epsilon,t} - \widetilde{\omega}_{\epsilon,t} \rvert  \leq c_4^*  t^{\alpha}.
\end{equation*}
And, for all $w_*=\omega(\xi), \xi \in  V_{\epsilon,t} \cap \widetilde{V}_{\epsilon,t}$, there exist $x_0, \widetilde{x}_0 \in V$, such that $\omega_0(x_0) = w_* = \widetilde{\omega}_0(\widetilde{x}_0)$, and we have
\begin{equation*}
       \lvert   x_0 - \widetilde{x}_0 \rvert \leq \frac{4 c_4^*  t^{\alpha}}{\Theta_1},
\end{equation*}
here $\lvert V \setminus (V_{\epsilon,t} \cap \widetilde{V}_{\epsilon,t}) \rvert \to 0 $ as $\gamma \to 0$. 

Thus, if time step  $t$ and disturbance parameter  $\epsilon$ are small enough, $\omega$ satisfies the Kolmogorov's non-degeneracy condition, 
the generating function representation of symplectic algorithm \eqref{2} has invariant tori on a set of large  measures. 
The generating function representation of phase flow in 
nearly integrable Hamiltonian system \eqref{xy3} also has invariant tori on a set of large  measures, and there is only a slight deformation in the common area. The Hausdorff's distance between the invariant tori in the phase spaces of the two systems is about $o(t^{\alpha})$.

\begin{remark}
   There are still some problems to be studied, such as the effect of $\Theta_1$ in Kolmogorov's non-degeneracy condition in Remark \ref{rem2} and $\beta$ in R\"{u}ssmann's non-degeneracy condition in Remark \ref{re1} on the whole system and invariant tori, which will be considered in the subsequent study.
\end{remark}

\begin{remark}
   Notice that not all sufficiently small time step sizes are appropriate. The structure of the set of time step sizes under the Kolmogorov's non-degeneracy condition is given by \cite{r7}. It turns out that the set of time step sizes is a Cantor-like set, and the density of the Cantor set is 1 at the origin of the real line. Although the frequency vector will drift under the R\"{u}ssmann's non-degeneracy condition, we can get very similar results after an analysis exactly like that in \cite{r7}. That is, the set is also a Cantor set and the density of which is also 1 at the origin of the real line. For details, please refer to \cite{r7}.
\end{remark}

%
%


\begin{thebibliography}{10}




\bibitem{r1}
Shang, Zaijiu.
\newblock {KAM} theorem of symplectic algorithms for {H}amiltonian systems.
\newblock {Numerische Mathematik}, 83(3):477--496, 1999.


\bibitem{r2}
Xu, Junxiang and You, Jiangong and Qiu, Qingjiu.
\newblock Invariant tori for nearly integrable {H}amiltonian systems with degeneracy.
\newblock {Mathematische Zeitschrift}, 226(3):375--387, 1997.

\bibitem{r4}
P{\"o}schel, J{\"u}rgen.
\newblock A lecture on the classical {KAM} theorem.
\newblock {Proceedings of Symposia in Pure Mathematics}, 69:707--732, 2001.

\bibitem{r5}
R{\"u}ssmann, Helmut.
\newblock Invariant tori in non-degenerate nearly integrable {H}amiltonian systems.
\newblock {Regular and Chaotic Dynamics}, 6(2):119--204, 2001.

\bibitem{r6}
Ding, Zhaodong and Shang, Zaijiu.
\newblock Numerical invariant tori of symplectic integrators for integrable {H}amiltonian systems.
\newblock {Science China Mathematics}, 61(9):1567--1588, 2018.

\bibitem{r7}
Shang, Zaijiu.
\newblock Resonant and {D}iophantine step sizes in computing invariant tori of {H}amiltonian systems.
\newblock {Nonlinearity}, 13(1):299--308, 2000.

\end{thebibliography}


\end{document}